\documentclass[a4paper,11pt]{article}
\usepackage[dvips]{graphicx,color}
\usepackage{amssymb}
\usepackage{amsthm}
\usepackage{amsmath}
\usepackage[ruled,vlined,linesnumbered]{algorithm2e}
\usepackage{index}
\usepackage{tikz}
\usepackage{pgfkeys}
\usepackage{pgfplots}
\usetikzlibrary{arrows,decorations.pathmorphing,backgrounds,positioning}
\usetikzlibrary{fit,petri,automata,chains}
\usepgfplotslibrary{statistics}

\DeclareMathOperator{\tr}{tr}

\begin{document}
\title{A Power Method for Computing Singular Value Decomposition}
\author{Doulaye Demb\'el\'e}
\date{Institut de G\'en\'etique et de Biologie Moleculaire et
  Cellulaire\\
  CNRS UMR7104, INSERM U964 and Universit\'e de Strasbourg\\
  67400 Illkirch, France~}
\maketitle

\begin{abstract}
The singular value decomposition (SVD) allows to write a 
matrix as a product of a left singular vectors matrix, a nonnegative
singular values diagonal matrix and a right singular vectors matrix. 
Among the applications of the SVD are the principal component analysis, 
the low-rank matrix approximation and the solving of a linear system of 
equations. The methods used for computing this decomposition
allow to get the complete or partial result.
For very large size matrix, the probabilistic methods 
allow to get partial result by using less computational load.
A power method is proposed in this paper for computing all or the $k$ first 
largest SVD subspaces for a real-valued matrix.
The $k$ first right singular vectors of this method are the $k$
columns of a neural network encoder weight matrix. The accuracy of
this iterative search method depends on the behavior of the singular
values and the settings of the gradient search optimizer used.
A R package implementing the proposed method is available at 
https://cran.r-project.org/web/packages/psvd/index.html.
\end{abstract}

\section{Introduction}
Nowadays, we are dealing with massive data in our current activities: words
count in text documents, images, gene expression from sequencing analyses, 
etc. Such data are typically organized in a matrix of size $m\times n$ where 
$m$ and $n$ are the number of rows and columns, respectively. 
We assume that $m\geq n$, this condition can be met by using a matrix 
transpose. A compact form of a large size data is often required for storing 
purpose or for revealing some links between the samples. The singular 
value decomposition (SVD) allows to express a real-valued matrix as a product 
of a left orthogonal, a diagonal and a right orthogonal matrices,
\cite{Eckart-al-1936,Stewart-1993}. The entries of the diagonal matrix are the
nonnegative-valued singular values and the columns of the orthogonal matrices are
the singular vectors. The properties of these vectors 
allow to write the matrix as a sum of rank one matrices.
The contribution of each term in this sum depends on the weight
of the associated singular value. 

One strategy for obtaining the SVD result consists of computing a symmetric 
matrix followed by some transformations. In a first category of 
transformation methods, a 
bidiagonalization or rank revealing decomposition is done followed by
a search for the factor matrices
\cite{Golub-al-1965,Golub-al-1970,Chan-1982,Demmel-al-1990,Demmel-al-1999}.
All singular values and singular vectors are computed at a time using 
this category of methods.
The quantification in SVD perturbation uncertainties is addressed 
in \cite{Yang-al-2024} while the convex set-oriented SVD method 
proposed in \cite{Fan-al-2024} overcomes the limitation of the nominal 
SVD method in the presence of uncertainties.
In a second category of transformation methods, an iterative search allows
to get the factor matrices entries
\cite{Lanczos-1950,Hestenes-al-1951,Arnoldi-1951,Saad-1980,
Knyazev-al-1994,Zhou-al-2007}.
This second category of methods allows to compute the $k$ largest or 
smallest singular values and singular vectors.
Another strategy for computing the SVD consists of a direct algebraic
manipulation of the matrix entries. The refinement method proposed 
in \cite{Ogita-al-2020} is based on matrices multiplication for computing 
the SVD of a full columns rank matrix.
The randomized methods allow a direct search for an 
orthogonal matrix associated with the first $k$ largest singular values,
\cite{Frieze-al-2004,Halko-al-2011,Gu-2015}. 
When the data matrix size is large, the first $k$ largest subspaces
are often enough to have a good data matrix approximation. 
However, to have accurate $k$ first singular values using the randomized 
method, a matrix of order greater $k$ should be used, 
see Table \ref{tab-svd-d3}. 

The neural network (NN) model is nowdays used in many applications,
\cite{Bishop-2006,LeCun-al-2015,Zhao-al-2024,Li-al-2025}. This 
model has a weight matrix which entries are computed from a training dataset.
The NN autoencoder model consists of two functions (encoding and decoding) 
allowing to store internally the input and then output it in a desirable format. 
This model has numerous applications: linear algebra, data 
denoising, regression analysis, texts and images generation, natural language 
processing, supervised and unsupervised classification. For the unsupervised 
classification or data clustering, the autoencoder model allows to get lower 
dimensional or denoised data which are used to obtain clusters. The singular 
values and the singular vectors of the SVD are the encoding results, while the 
decoding can consist of using a reconstruction form the $k$ first singular 
values and singular vectors.
A deep embedding clustering (DEC) method is proposed in \cite{Xie-al-2016}
by using an autoencoder to have a lower dimensional dataset, the 
decoder output. The cluster centers in DEC are learned from this lower 
dimensional data after a minimization of an objection function involving a 
Kullback-Leibler divergence and the t-distribution. 
The deep divergence clustering method proposed in \cite{Kampffmeyer-al-2019} 
uses convolutional layers, autoencoder, for images datasets. An autoencoder is 
used in the deep-k-means method proposed in \cite{Fard-al-2020}. A dynamic 
autoencoder is used in \cite{Mrabah-al-2020} for clustering a dataset by gradually 
and smoothly eliminating the reconstruction objective function in favor of a 
construction one. The autoencoder model parameters and the clustering 
parameters are simultaneously learned in the method proposed in 
\cite{Boubekki-al-2021} by using a certain class of gaussian mixtures based 
objective function.  The deep convolutional embedded clustering algorithm 
proposed in \cite{Lu-al-2022} uses a convolutional 
autoencoder followed by convolutional neural network. The unsupervised 
companion objectives method proposed in \cite{Trosten-al-2024} allows 
to reduce the mismatch which occurs when the optimization of an objective 
(the autoencoder reconstruction loss function) has a negative impact on the 
optimization of another objective (the clustering loss function). 

In data matrix approximation applications, it has been shown that the singular 
vectors can be the NN model weight matrix columns, 
\cite{Oja-1982,Bourlard-al-1988,Baldi-al-1989,Kramer-1991,Xu-1993}.
An independent optimization method is very often used to get the singular vectors.
For the iterative search method proposed, the 
optimizer output singular vectors are the NN weight matrix columns.
The right singular vectors matrix of the probabilistic method correspond to
an NN weight matrix. However, the probabilistic method results are less accurate
when the number of the singular values increases.
The proposed algorithm is inspired from the probabilistic 
random search approaches, \cite{Halko-al-2011,Gu-2015}. 
The difference comes from a minimization of an objective function and 
the use of a power method to get simultaneously the singular vectors. 
The proposed method is more accurate, can be used for any matrix and also 
allows to compute the complete SVD result.

\subsection{Computing the SVD}
The SVD of any matrix $\mathbf{X}\in\mathbb{R}^{m\times n}$ is, 
\cite{Marcus-al-1964,Golub-al-1996,Horn-al-2019}:
\begin{equation}
  \mathbf{X}=\mathbf{UD}\mathbf{W}^{\mathsf{T}} \label{eq-svd}
\end{equation}
where $\mathbf{U}$ and $\mathbf{W}$ are the left and right
orthogonal matrices of
order $m$ and $n$, respectively, $\mathbf{D}$ is a matrix of size $m\times n$.
The diagonal entries of the upper square part of the matrix 
$\mathbf{D}$ are nonnegative, the singular values:
$d_{11}\geq d_{22}\geq \ldots \geq d_{nn}$. In the thin SVD, the matrix
$\mathbf{U}$ is of size $m\times n$ and $\mathbf{D}$ is
of order $n$. Let us observe that $-\mathbf{U}$ and $-\mathbf{W}$ are jointly
also solution of the singular value decomposition in relation (\ref{eq-svd}).
The column orthogonal property of the matrices $\mathbf{U}$ and
$\mathbf{W}$ allows to write:
$\mathbf{U}^{\mathsf{T}}\mathbf{U}=\mathbf{W}^{\mathsf{T}}\mathbf{W}=\mathbf{I}_n$,
where $\mathbf{I}_n$ is a $n$-order identity matrix.
Using the orthogonal property, relation (\ref{eq-svd}) can be written as:
\begin{equation}
  \mathbf{X}=\sum_{j=1}^{n}d_{jj}\mathbf{u}_j\mathbf{w}_j^{\mathsf{T}}
  \label{eq-svd2}
\end{equation}
where the unit $\ell_2$-norm left and right singular
vectors $\mathbf{u}_j$ and $\mathbf{w}_j$ are the $j$-th
columns of the matrices $\mathbf{U}$ and $\mathbf{W}$, respectively,
$\mathbf{u}_k^{\mathsf{T}}\mathbf{u}_j = 1$, if $k=j$, $0$
otherwise, idem for vectors $\mathbf{w}_j$. Again, for column $j$, the
couples of vectors ($\mathbf{u}_j, \mathbf{w}_j$) and ($-\mathbf{u}_j,
-\mathbf{w}_j$) have the same contribution. If there are $n$ nonzero
singular values $d_{jj}$, then the matrix $\mathbf{X}$ is of full
rank. Otherwise, the number of nonzero scalars $d_{jj}$ is
the actual rank of $\mathbf{X}$. The sum in relation
(\ref{eq-svd2}) is often truncated to $r$ terms ($r<n$) leading to a
lower rank matrix approximation. The percentage of
the information associated with the $r$ first largest singular values is:
\begin{equation}
reconst.(\%) =
100\sum_{j=1}^rd_{jj}\left/\sum_{j=1}^nd_{jj}\right.
\label{eq-svd-reconstr}
\end{equation}
where the parameter $r$ is chosen to have the desired reconstruction rate.

\subsection{Eigenvalue decomposition}
The eigenvalue decomposition (EVD) or eigendecomposition 
can be performed for any square matrix, \cite{Golub-al-2000} .
The EVD of the symmetric matrix $\mathbf{X}^{\mathsf{T}}\mathbf{X}$ is:
\begin{equation}
  \mathbf{X}^{\mathsf{T}}\mathbf{X} =
  \mathbf{V}\boldsymbol{\Lambda}\mathbf{V}^{\mathsf{T}}  \label{eq-eigen}
\end{equation}
where $\boldsymbol{\Lambda}$ is a diagonal matrix which entries are 
nonnegative real-valued scalars, the eigenvalues:
$\lambda_1\geq\lambda_2\geq\ldots\geq\lambda_n\geq 0$,
$\mathbf{V}$ is an orthogonal matrix,
$\mathbf{V}^{\mathsf{T}}\mathbf{V}=\mathbf{I}_n$.
Like in relation (\ref{eq-svd2}), the EVD can be
expressed as the sum of rank one matrices, the spectral decomposition:
\begin{equation}
  \mathbf{X}^{\mathsf{T}}\mathbf{X} =
  \sum_{j=1}^n\lambda_j\mathbf{v}_j\mathbf{v}_j^{\mathsf{T}}\label{eq-eigen2a}
\end{equation}
where $\mathbf{v}_j$ is the $j$-th column of $\mathbf{V}$, a unit $\ell_2$-norm
vector. The above sum can also be truncated to get a $r<n$ low-rank
matrix $\mathbf{X}^{\mathsf{T}}\mathbf{X}$ approximation.
In the principal component analysis, see below, the columns of the matrix
$\mathbf{V}$ correspond to the projection subspaces. One can observe that
$-\mathbf{V}$ is also solution of the EVD, relation (\ref{eq-eigen}).
This observation allows to fulfil the Perron-Frobenius theorem,
\cite{Perron-1907,Frobenius-1912,Meyer-2000}, which stipulates that the entries of
the main eigenvector of a nonnegative-valued matrix are all positive or zero.
Nonnegative-valued matrices occur for many real-life data, e.g. genes
expression and images.
The main eigenvector is used in the ranking applications,
\cite{Brin-al-1998,Langville-al-2005}.
To compute the eigenvalue with maximum modulus and the associated eigenvector,
the iterative power method, \cite{Golub-al-1996}, is very often used, see
algorithm \ref{algo-power} in the Annexes.
The convergence rate of the power method depends on the modulus of the second
eigenvalue compared to the first, \cite{Golub-al-1996,Dembele-2021}.

\subsection{Singular value decomposition and eigenvalue decomposition}
It has been shown that the SVD factor matrices of a symmetric matrix
are related to those of the EVD 
\cite{Golub-al-1965,Hanson-al-1969,Stewart-1993}. Indeed,
using the expression (\ref{eq-svd}) of $\mathbf{X}$ in the relation
(\ref{eq-eigen}) allows to write:
\begin{equation}
  \mathbf{X}^{\mathsf{T}}\mathbf{X} = (\mathbf{UDW}^{\mathsf{T}})^{\mathsf{T}}
  \mathbf{UDW}^{\mathsf{T}} = \mathbf{WD}^2\mathbf{W}^{\mathsf{T}}\nonumber
\end{equation}
Comparing the above expression to relation (\ref{eq-eigen}),
the SVD factors using the EVD factors are:
\begin{equation}
  \mathbf{W}=\mathbf{V}\mbox{ ; }
  \mathbf{D}=\boldsymbol{\Lambda}^{1/2}\mbox{ and }
\mathbf{U} = \mathbf{XW}\boldsymbol{\Lambda}^{-1/2}\label{eq-svd-eigen}
\end{equation}
The right singular vectors of $\mathbf{X}$ are the eigenvectors of 
$\mathbf{X}^{\mathsf{T}}\mathbf{X}$. The entries of the diagonal matrix 
$\boldsymbol{\Lambda}$, the eigenvalues, are the diagonal
entries of the symmetric matrix $(\mathbf{XW})^{\mathsf{T}}\mathbf{XW}$.

\subsection{Bidiagonalization methods for computing the SVD}
The SVD is computed in \cite{Golub-al-1996} by using a bidiagonalization
of the matriox $\mathbf{X}$ or a QR decomposition of the matrix
$\mathbf{X}^{\mathsf{T}}\mathbf{X}$. In the first case, a sequence of
Householder transformations leads to an upper bidiagonal matrix 
$\mathbf{B}$ and two other matrices, $\mathbf{U}_1$ and $\mathbf{V}_1$, 
\cite[page 252]{Golub-al-1996}. In the second case,
the symmetric matrix allows to get a column orthogonal matrix $\mathbf{Q}$ 
and an upper triangular matrix $\mathbf{R}$, \cite[page 448]{Golub-al-1996}. 
In both cases, another manipulation of the factor matrices is used to get 
the SVD. The lower triangular, diagonal and upper
triangular (LDU) decomposition is used in
\cite{Demmel-al-1999} to obtain the SVD. The diagonal entries of the
triangular matrices $\mathbf{L}$ and $\mathbf{U}$ are all $1$. The
diagonalization of the matrix $\mathbf{X}$ before obtaining the SVD
leads to relatively accurate singular values. All singular values and
the singular vectors are simultaneously recovered by these approaches.

\subsection{Iterative methods for computing the SVD}
The iterative search methods are used to construct a collection of 
mutually orthogonal subspaces leading to the singular values and 
singular vectors. From the Rayleigh-Ritz theorem, the minimum and the
maximum eigenvalue of a symmetric matrix $\mathbf{A}$
verify ($\mathbf{v}^{\mathsf{T}}\mathbf{v}=1$): 
$\lambda_{max}(\mathbf{A}) \geq \mathbf{v}^{\mathsf{T}}\mathbf{A}\mathbf{v}
\geq \lambda_{min}(\mathbf{A})$.
For a non null vector $\mathbf{v}$, the Rayleigh quotient is: 
$r(\mathbf{v}) = \mathbf{v}^{\mathsf{T}}\mathbf{A}\mathbf{v}/
\mathbf{v}^{\mathsf{T}}\mathbf{v}$. A gradient search method based on
the Rayleigh quotient has been used in \cite{Hestenes-al-1951} to get
the eigenvalues and the eigenvectors of the matrix $\mathbf{A}$. 

The Krylov $k$ columns (subspaces) matrix associated with the symmetric
matrix $\mathbf{A}$ and the vector $\mathbf{v}$ is: 
$\mathsf{K}(\mathbf{A}, \mathbf{v},k) = (\mathbf{v},\ \mathbf{Av},\
\mathbf{A}^2\mathbf{v},\ \ldots,\ \mathbf{A}^{k-1}\mathbf{v})$. A next
subspace, a column, in this matrix is obtained from the previous by using 
only a matrix-vector product. The Lanczos iterative methods consist of 
generating Krylov subspaces to get the eigenvalues and the eigenvectors. 
The Arnoldi method uses minimum iterations 
\cite{Arnoldi-1951,Knyazev-al-1994}. 
The SVD is also used to solve the linear system of equations which
occurs in many engineering problems: $\mathbf{Ax} = \mathbf{b}$, where
$\mathbf{A}$ and $\mathbf{b}$ are known matrix and vector of size
$m\times n$ and $n$, respectively, and $\mathbf{x}$ is the unknown vector.
For a square matrix, $\mathbf{A}$, the
Lanczos algorithm, \cite{Lanczos-1950}, allows to get a solution
through the EVD. For a non-square matrix $\mathbf{A}$, the normal equation,
$\mathbf{A}^{\mathsf{T}}\mathbf{Ax} = \mathbf{A}^{\mathsf{T}}\mathbf{b}$,
involves a symmetric matrix which can be used to search for
$\mathbf{x}$. For a sparse matrix with large size, the
Lanczos algorithm has been extended by Larsen, \cite{Larsen-1998}, to
bidiagonalization of the non-square matrix leading to perform its SVD.
For high dimensional data, when the purpose is
to get a low-rank data matrix approximation, it may not be necessary to
get all singular values. The Lanczos algorithm extension proposed in
\cite{Baglama-al-2005} allows to compute the few first largest or
smallest singular values and the singular vectors.
In the Davidson type iterative methods, polynomial filters are
used to get the SVD subspaces, \cite{Sleijpen-al-1996,Zhou-al-2007}.

\subsection{Probabilistic SVD computation}
The computation of the SVD for massive data has been considered in
\cite{Liberty-al-2007,Halko-al-2011,Gu-2015,Tropp-al-2018}, where,
only the $r$ first largest singular values and singular vectors are
computed using a probabilistic (randomized) approach. The search
consists of finding column orthogonal matrix $\mathbf{Q}$ of size
$m\times r$ such that:  
\begin{equation}
\|\mathbf{X} - \mathbf{QQ}^{\mathsf{T}}\mathbf{X}\| \leq \epsilon
\label{eq-proba-svd-a}
\end{equation}
where $\|.\|$ is the $\ell_2$-norm, $r$ and $\epsilon$
are settings for the rank and the tolerance of the approximation.
In the randomized approach, the search for the $r$ first singular values 
and singular vectors are performed through a random matrix
$\boldsymbol{\Omega}$ of size, $n\times (r+\ell)$,
where $\ell$ is an integer chosen such that $r+\ell<n$.
The entries of this matrix are typically from the standard
normal distribution: $\boldsymbol{\Omega}=(\omega_{ij})$,
$\omega_{ij}\sim \mathcal{N}(0,1)$, $\forall i,j$.
$\boldsymbol{\Omega}$ is used to define an intermediate matrix,
$\mathbf{Y} \triangleq \mathbf{X}\boldsymbol{\Omega}$, which has few columns
compared to $\mathbf{X}$. A QR decomposition, 
\cite{Golub-al-1996}, of the matrix $\mathbf{Y}$ is next performed:
$\mathbf{X}\boldsymbol{\Omega} = \mathbf{QR}$, where the columns of the
matrix $\mathbf{Q}$ are orthogonal.
Setting $\mathbf{Z}\triangleq \mathbf{Q}^{\mathsf{T}}\mathbf{X}$ allows to 
have $\mathbf{QZ} = \mathbf{QQ}^{\mathsf{T}}\mathbf{X}$ or
$\mathbf{QZ}=\mathbf{X}$, thank to the column orthogonality of $\mathbf{Q}$.
Thus, the SVD of the matrix
$\mathbf{Z}$ is performed and the result $\mathbf{QZ}$ is returned.
Indeed, $\mathbf{QZ}=\mathbf{Q}\tilde{\mathbf{U}}\mathbf{DW}^{\mathsf{T}}
=\mathbf{UDW}^{\mathsf{T}}$.
To have rapid singular values decay, the intermediate matrix is computed as:
$\mathbf{Y}\triangleq (\mathbf{XX}^{\mathsf{T}})^q\mathbf{X}\boldsymbol{\Omega}$.
Indeed, using the relation $\mathbf{X}=\mathbf{UDW}^{\mathsf{T}}$,
the column orthogonal property of the
matrices $\mathbf{U}$ and $\mathbf{W}$, and for an integer $q\geq 0$,
allow to write:
\begin{equation}
  (\mathbf{XX}^{\mathsf{T}})^q\mathbf{X} =
  \mathbf{U}\mathbf{D}^{2q+1}\mathbf{W}^{\mathsf{T}}  \label{eq-proba-svd}
\end{equation}
This relation shows a ($2q+1$)-power of the singular values for the matrix
$\mathbf{Y}$.
The rapid singular values decay leads to a fast convergence algorithm,
\cite{Halko-al-2011,Gu-2015}. A historical overview of the randomized 
computing and its applications in linear algebra problems is provided in
\cite{Halko-al-2011}.

Using the probabilistic SVD method, the $r$ rank approximation error
is slightly greater than the bound given below in relation
(\ref{eq-rank-r-approx}), \cite{Halko-al-2011}. In addition, the number of columns,
$r+\ell$, of the random matrix $\boldsymbol{\Omega}$ should be large if
accurate values are expected for the first $r$ singular values.
The proposed method leads to more accurate $r$ first singular values by
using a $r$ columns random matrix. More over, the proposed method
can be used to compute the complete SVD even for a small size matrix.

\section{Computing the SVD by using a power method}
Thank to relation (\ref{eq-eigen2a}), the power method can be used to get
all eigenvalues and eigenvectors, \cite{Golub-al-1996}. Indeed, let us define
$\mathbf{S}_1 = \mathbf{X}^{\mathsf{T}}\mathbf{X}$. When the first eigenvector 
is obtained, a new matrix can be formed as
$\mathbf{S}_2 = \mathbf{S}_1-\lambda_1\mathbf{v}_1\mathbf{v}_1^{\mathsf{T}}$.
The matrix $\mathbf{S}_2$ is used to get the second eigenvector and associated
eigenvalue. This process is repeated to get all the eigenvector and eigenvalue
pairs. That is the Jordan SVD algorithm, \cite{Stewart-1993}.
Instead of computing the eigenpairs of the matrix 
$\mathbf{X}^{\mathsf{T}}\mathbf{X}$, one by one, the proposed algorithm
computed them simultaneously. The column orthogonal property
of the matrix $\mathbf{V}$ allows to have: 
\begin{equation}
(\mathbf{X}^{\mathsf{T}}\mathbf{X})^q =
  \mathbf{V}\boldsymbol{\Lambda}^q\mathbf{V}^{\mathsf{T}}\label{eq-eigen3}
\end{equation}
where $q$ is an integer, $q\geq 1$. That means, all power $q$ of the matrix
$\mathbf{X}^{\mathsf{T}}\mathbf{X}$ have the same eigenvectors, only the
eigenvalues differ by the power $q$ and we also have:
$\lambda_1^q\geq\lambda_2^q\geq\ldots\geq\lambda_n^q\geq 0$.
At least, the first eigenvalue is nonzero, that allows to write:
\begin{eqnarray}
  (a)&:&1\geq (\lambda_2/\lambda_1)\geq (\lambda_3/\lambda_1)\geq
  \ldots\geq(\lambda_n/\lambda_1)\\
  (b)&:&1\geq (\lambda_2/\lambda_1)^q\geq (\lambda_3/\lambda_1)^q\geq
  \ldots\geq(\lambda_n/\lambda_1)^q \\
  (c)&:&(\lambda_j/\lambda_1)^q\leq (\lambda_j/\lambda_1)\mbox{ ; }
  j=1,2,\ldots,n \label{eq-decay}
\end{eqnarray}
The proposed power method is based on the 
Frobenius norm of the matrix $\mathbf{X}$:
\begin{equation}
  \|\mathbf{X}\|_{\mathsf{F}}^2 = \sum_{i=1}^m\sum_{j=1}^nx_{ij}^2
  = \tr\left[\mathbf{X}^{\mathsf{T}}\mathbf{X}\right]
  = \tr\left[\mathbf{XX}^{\mathsf{T}}\right]
  = \sum_{j=1}^nd_{jj}^2\label{eq-Frobenius-norm}
\end{equation}
where $\tr[.]$ is the matrix trace operator. Let $\mathbf{X}_r$ be the
matrix obtained using the $r$ first terms in the sum of relation
(\ref{eq-svd2}). The Frobenius norm allows to quantify the error of
the $r$-rank $\mathbf{X}$ approximation,
\cite{Golub-al-1996,Horn-al-2019}: 
\begin{equation}
  \|\mathbf{X} - \mathbf{X}_r\|_{\mathsf{F}}^2 = \sum_{j=r+1}^nd_{jj}^2
  \label{eq-rank-r-approx}
\end{equation}

Relation (\ref{eq-svd-eigen}) shows a one-to-one relationship between the
SVD factors of $\mathbf{X}$ and the EVD factors
of $\mathbf{X}^{\mathsf{T}}\mathbf{X}$. From relation (\ref{eq-svd-eigen}),
the SVD of the matrix $\mathbf{X}$ can be written as:
\begin{equation}
  \mathbf{X} = \mathbf{XWW}^{\mathsf{T}}\label{eq-x-xwwt}
\end{equation}
When $\mathbf{W}\triangleq \mathbf{Q}$, the above relation is
equivalent to that is used in relation (\ref{eq-proba-svd-a}) and is
obtained when performing the EVD of the matrix $\mathbf{XX}^{\mathsf{T}}$.
Since we assumed that $n\leq m$, working with the matrix
$\mathbf{X}^{\mathsf{T}}\mathbf{X}$ necessitates less computational load
compared to the matrix $\mathbf{XX}^{\mathsf{T}}$.
To obtain the column orthogonal matrix $\mathbf{W}$ given
$\mathbf{X}$, a norm of the difference between the left-hand and
right-hand side terms in the above relation can be minimized. In the
probabilistic SVD approach, a column orthogonal matrix of size
$r+\ell$, $\ell>r$, is used and the QR decomposition is performed
followed by matrices multiplication to get the SVD. Here, the
entries of a random column orthogonal matrix of size $n\times r$ are
iteratively modified to fit the factor $\mathbf{W}$ of the SVD.
The algorithm used is similar to the power
method but the $r\leq n$ eigenvectors are simultaneously recovered. 

\subsection{Search for the projection matrix by using the gradient method}
To obtain the matrix $\mathbf{W}$ in relation (\ref{eq-x-xwwt}), the
Frobenius norm of the difference between $\mathbf{X}$ and
$\mathbf{XWW}^{\mathsf{T}}$ is minimized, subject to the column orthogonality
constraint for the matrix $\mathbf{W}$. 
\begin{equation}
\min\mathsf{J}(\mathbf{W})\mbox{; }
\mathsf{J}(\mathbf{W}) = \frac{1}{2}\left\|\mathbf{X} - \mathbf{XWW}^{\mathsf{T}}
\right\|_{\mathsf{F}}^2 \mbox{ subject to: }
\mathbf{W}^{\mathsf{T}}\mathbf{W}=\mathbf{I}\label{eq-optim}
\end{equation}
The Lagrangian of this optimization problem is:
\begin{equation}
  \mathcal{L}(\mathbf{W},\Upsilon) =
  \frac{1}{2}\left\|\mathbf{X}-\mathbf{XWW}^{\mathsf{T}}\right\|_{\mathsf{F}}^2-
  \left\|\Upsilon\left(\mathbf{W}^{\mathsf{T}}\mathbf{W}-
  \mathbf{I}\right)\right\|_{\mathsf{F}}^2\label{eq-lagrangian}
\end{equation}
where the entries of the matrix $\Upsilon$ are the Lagrange multipliers
associated with the equality constraints.
To minimize $\mathsf{J}(\mathbf{W})$, the gradient descent method is
used, the solution at iteration $t$ is obtained after
updating the previous one: $\mathbf{W}^{(t)} = \mathbf{W}^{(t-1)} -
\eta\nabla\mathcal{L}(\mathbf{W},\Upsilon)$, where $\eta$ is a
positive step parameter and $\nabla\mathcal{L}(\mathbf{W},\Upsilon)$
is the gradient of the Lagrangian. The derivative of
$\mathcal{L}(\mathbf{W},\Upsilon)$ with respect to $\mathbf{W}$ and
$\Upsilon$ are then required. 

Relation (\ref{eq-Frobenius-norm}) allows to write:
\begin{eqnarray}
  \left\|\mathbf{X} - \mathbf{XWW}^{\mathsf{T}}\right\|_{\mathsf{F}}^2 &=&
  \tr\left[(\mathbf{X}-\mathbf{XWW}^{\mathsf{T}})(\mathbf{X}-
    \mathbf{XWW}^{\mathsf{T}})^{\mathsf{T}}\right]\nonumber\\
  &=& \tr\left[\mathbf{X}\mathbf{X}^{\mathsf{T}} -
    2\mathbf{XWW}^{\mathsf{T}}\mathbf{X}^{\mathsf{T}}
    + \mathbf{XWW}^{\mathsf{T}}\mathbf{WW}^{\mathsf{T}}\mathbf{X}^{\mathsf{T}}
    \right]\nonumber
\end{eqnarray}
For any column orthogonal matrix $\mathbf{W}$ of size $n\times r$,
$\mathbf{W}^{\mathsf{T}}\mathbf{W}=\mathbf{I}_r$, that allows to write:
\begin{eqnarray}
  \left\|\mathbf{X} - \mathbf{XWW}^{\mathsf{T}}\right\|_{\mathsf{F}}^2
  &=& \tr\left[\mathbf{X}\mathbf{X}^{\mathsf{T}} 
    -\mathbf{XWW}^{\mathsf{T}}\mathbf{X}^{\mathsf{T}} \right]\label{eq-frob-tr}\\
  \left\|\mathbf{W}^{\mathsf{T}}\mathbf{W}-\mathbf{I}\right\|_{\mathsf{F}}^2
  &=& 0\label{eq-frob-tr2}
\end{eqnarray}
For any $n$-order orthogonal matrix $\mathbf{W}$, we have
$\mathbf{W}\mathbf{W}^{\mathsf{T}} = \mathbf{I}_n$. However, all
$n$-order orthogonal matrices do not allow to verify relation
(\ref{eq-eigen}) for a given matrix $\mathbf{X}$. In addition, we
would like to be able to compute a column orthogonal matrix
$\mathbf{W}$ of size $n\times r$, $r<n$. In this case, we still have
$\mathbf{W}^{\mathsf{T}}\mathbf{W} = \mathbf{I}_r$ but
$\mathbf{W}\mathbf{W}^{\mathsf{T}} \neq \mathbf{I}_n$. The matrix
$\mathbf{XX}^{\mathsf{T}}-\mathbf{XWW}^{\mathsf{T}}\mathbf{X}^{\mathsf{T}}$
is symmetric and has nonnegative-valued
eigenvalues. The trace of this matrix is the sum of its
eigenvalues. For a $n$-order column orthogonal matrix $\mathbf{W}$,
relation (\ref{eq-frob-tr}) shows that the global minimum of the objective
function $\mathsf{J}(\mathbf{W})$ is zero. For a $r$-column
orthogonal matrix $\mathbf{W}$, $r<n$, we have:
\[tr\left[\mathbf{XX}^{\mathsf{T}} -
    \mathbf{XWW}^{\mathsf{T}}\mathbf{X}^{\mathsf{T}}\right] =
  tr\left[\mathbf{XX}^{\mathsf{T}}\right] -
  tr\left[\mathbf{XW}(\mathbf{XW})^{\mathsf{T}}\right] =  \sum_{j=r+1}^nd_{jj}\]
where $\mathbf{XW}(\mathbf{XW})^{\mathsf{T}}$ is
the diagonal matrix which entries are the $r$ first singular values,
see the text under relation (\ref{eq-svd-eigen}).
The relation (118) in
\cite[page 13]{Petersen-al-2012}, can be written as:
\[\frac{\partial}{\partial\mathbf{W}}
\tr\left[\mathbf{AWBW}^{\mathsf{T}}\mathbf{C}\right] =
\mathbf{A}^{\mathsf{T}}\mathbf{C}^{\mathsf{T}}\mathbf{WB}^{\mathsf{T}} +
\mathbf{CAWB}\]
Using this relation with $\mathbf{A}=\mathbf{X}$,
$\mathbf{B}=\mathbf{I}$ and $\mathbf{C}=\mathbf{X}^{\mathsf{T}}$, the
gradient of the objective function $\mathsf{J}(\mathbf{W})$ is:
\begin{equation}
  \frac{\partial\mathsf{J}(\mathbf{W})}{\partial\mathbf{W}} =
  -\mathbf{X}^{\mathsf{T}}\mathbf{XW}\nonumber
\end{equation}
The matrix $\mathbf{W}$ in the above relation is that available in the
previous iteration, $\mathbf{W}^{(t-1)}$. Then, the gradient updating
equation is: $(\mathbf{I}_n +
\eta\mathbf{X}^{\mathsf{T}}\mathbf{X})\mathbf{W}^{(t-1)}$. Therefore, 
at iteration $t=1,2,\ldots$ and given column orthogonal matrix
$\mathbf{W}^{(t-1)}$, two alternate steps are used to get the update
$\mathbf{W}^{(t)}$:
\begin{eqnarray}
\tilde{\mathbf{W}} &=& \mathbf{G}\mathbf{W}^{(t-1)}\label{eq-grad-a}\\
\mathbf{W}^{(t)} &=& \mbox{Gram-Schmidt}(\tilde{\mathbf{W}})  \label{eq-grad}
\end{eqnarray}
where
$\mathbf{G}\triangleq \mathbf{I}_n+\eta(\mathbf{X}^{\mathsf{T}}\mathbf{X})^q$,
$\eta$ is a positive step parameter and $q\geq 1$ is an integer.
The settings for $\eta$ and $q$ are discussed later.
The method to use for obtaining the column orthogonal matrix $\mathbf{W}$ is
summarized in algorithm \ref{algo-svd-gradient}.
\begin{algorithm}
\caption{Gradient search for column orthogonal matrix $\mathbf{W}$}
\label{algo-svd-gradient}
\SetAlgoLined
\KwData{matrix $\mathbf{X}$ of size $m\times n$, settings
$\eta$, $q$, $\epsilon$ and maximum number of iterations ($itmax$);}
\KwResult{column orthogonal matrix $\mathbf{W}$;}
Initialization: select randomly entries of a matrix $\mathbf{W}$ from a
  normal or an uniform distribution and perform orthogonalization
  $\mathbf{W}^{(0)}=$ Gram-Schmidt($\mathbf{W}$), compute
  $\mathbf{G} = \mathbf{I}_n+\eta(\mathbf{X}^{\mathsf{T}}\mathbf{X})^q$\;
\For{$t\leftarrow 1,2,\ldots, itmax$}{
  compute $\tilde{\mathbf{W}}$ using relation (\ref{eq-grad-a})\;
  compute $\mathbf{W}^{(t)}$ using relation (\ref{eq-grad})\;
  $\delta=\|\mathbf{W}^{(t)}-\mathbf{W}^{(t-1)}\|_{\mathsf{F}}^2$\;
  \If{$\delta\leq \epsilon$} {stop\;}
}
\end{algorithm}

The columns of the result matrix $\mathbf{W}$ are the eigenvectors matrix 
of the matrix $\mathbf{X}^{\mathsf{T}}\mathbf{X}$ and the right 
singular vectors matrix of $\mathbf{X}$.
Using relation (\ref{eq-eigen}), the eigenvalues of the matrix
$\mathbf{X}^{\mathsf{T}}\mathbf{X}$ are the entries of the diagonal matrix
$\boldsymbol{\Lambda}$ or the diagonal entries
of the symmetric matrix $(\mathbf{XW})^{\mathsf{T}}\mathbf{XW}$. 
Relation (\ref{eq-svd-eigen}) allows to get the SVD factor matrices 
$\mathbf{D}$ and $\mathbf{U}$.

\subsection{Relation with the power method}
Algorithm \ref{algo-svd-gradient} is a power method allowing to get
the $r\leq n$ eigenvectors simultaneously and not one by one like by the
classical power method in algorithm \ref{algo-power}, see Annexes. The
matrix $\mathbf{I}_n + \eta\mathbf{X}^{\mathsf{T}}\mathbf{X}$
is used in algorithm \ref{algo-svd-gradient} instead of the matrix
$\mathbf{X}^{\mathsf{T}}\mathbf{X}$. Indeed, 
the eigenvector $\mathbf{z}$ associated with the eigenvalue $\gamma$ of the
matrix $\mathbf{I}_n+\eta\mathbf{X}^{\mathsf{T}}\mathbf{X}$
allows to write:
\[\left(\mathbf{I}_n + \eta\mathbf{X}^{\mathsf{T}}\mathbf{X}
\right)\mathbf{z} = \gamma\mathbf{z} 
\Rightarrow \mathbf{X}^{\mathsf{T}}\mathbf{X}\mathbf{z}=
\left(\frac{\gamma-1}{\eta}\right)\mathbf{z}\]
The above relation means the matrices
$\mathbf{I}_n + \eta\mathbf{X}^{\mathsf{T}}\mathbf{X}$ and
$\mathbf{X}^{\mathsf{T}}\mathbf{X}$ have the same eigenvectors which are
associated with related eigenvalues, $j=1,2,\ldots,n$:
\begin{eqnarray}
  \lambda_j  &=& (\gamma_j-1)/\eta\label{eq-eigenGa}\\
  \gamma_j &=& \eta\lambda_j+1\label{eq-eigenGb}
\end{eqnarray}
where $\lambda_j$ and $\gamma_j$ are the eigenvalues of the
matrix $\mathbf{X}^{\mathsf{T}}\mathbf{X}$ and
$\mathbf{I}_n + \eta\mathbf{X}^{\mathsf{T}}\mathbf{X}$, respectively.
Relation (\ref{eq-eigenGb}) shows that the smallest eigenvalue of the matrix
$\mathbf{I}_n + \eta\mathbf{X}^{\mathsf{T}}\mathbf{X}$ is one, meaning
this matrix is always of full rank even for a rank deficient matrix
$\mathbf{X}$. A safe computation of the norm
$\|\mathbf{W}^{(t)} - \mathbf{W}^{(t-1)} \|_{\mathsf{F}}^2$, step 5 of
alogithm \ref{algo-svd-gradient}, at any iteration
requires to have the same number of columns for the matrix $\mathbf{W}$.
When the matrix $\mathbf{X}$ is not of full rank, there exit a matrix
$\mathbf{W}$ of size $n\times r$, $r<n$, allowing to have an exact
reconstruction of the $\mathbf{X}$. Therefore, using alone the matrix
$\mathbf{X}^{\mathsf{T}}\mathbf{X}$ in algorithm \ref{algo-svd-gradient} may
lead to a problem for computing the above norm if we are not careful. 

\subsection{Convergence analysis and choice of the settings $\eta$ and $q$}
\label{ssect-eta-q}
From the convex property of the objective function, the algorithm
\ref{algo-svd-gradient} is convergent. The rate of convergence depends
on the behavior of the singular values of the matrix $\mathbf{X}$, the 
precision $\epsilon$ and the step parameter $\eta$.
The rate of convergence for the probabilistic SVD is analyzed in
\cite{Halko-al-2011,Gu-2015} where faster convergence is observed for matrix
which singular values decay quickly. That means, larger the difference between
successive singular values is, faster is the convergence. This result
corresponds to the modulus conditions for the power method. Relation
(\ref{eq-decay}) shows a possibility to modify the singular values decay 
by using a power $q$ of the matrix in the EVD.
From relation (\ref{eq-eigenGb}), a small $\eta$ value (near zero)
may lead to near one eigenvalues for the matrix $\mathbf{G}$.
On the other hand, a value greater to $1$ for the parameter $\eta$ 
leads to increasing the eigenvalues of the matrix $\mathbf{G}$. 
Therefore, the choices $q>1$ and $\eta\geq 1$, likely lead to rapid
singular values decay for the matrix $\mathbf{G}$. However, when the
matrix $\mathbf{X}$ has large eigenvalues, small $\eta$ value
may be good.

The matrix 
$\mathbf{G}=\mathbf{I}_n+\eta(\mathbf{X}^{\mathsf{T}}\mathbf{X})^q$
which is computed once before starting the iterative search, may 
lead to numerical problem when $\mathbf{X}$ has some very 
small or high entries. At each iteration, a matrix multiplication and a
Gram-Schmidt (GS) orthogonalization are performed. The modified
GS orthogonalization method is summarized in algorithm
\ref{algo-Gram-Schmidt}, see Annexes.
The matrices multiplication and the GS orthogonalization are the
heavy computing time part of the proposed method which is implemented
in the \textit{psvd} R, \cite{Rproject}, package
(https://cran.r-project.org/web/packages/psvd/index.html). A C
code is used to speed-up the matrices multiplication and the GS
orthogonalization.

\section{Experimental results}
The implementation of algorithm \ref{algo-svd-gradient} in 
the R package \textit{psvd} has the 
following inputs:  $X$ the matrix, $r$ (the first singular values, 
default is 2), $\eta$ (the gradient search positive scalar, default 
is 1), $q$ (power integer, default is 2), \textit{itmax} (maximum 
iterations number, default is 200), $err$ (convergence error default 
is 1e-8) and \textit{mySeed} (a seed used for generating an initial 
solution, default is 50).  The matrix $X$ is mandatory and $r$ is a 
user supply parameter. 
The choice of the parameters $\eta$ and $q$ are discussed in 
subsection \ref{ssect-eta-q}. A large number of iterations 
number \textit{itmax} may lead to accurate singular values, 
especially a large size data matrix where the convergence can be 
slow or for a challenging matrix example. The parameter 
\textit{err} is used for convergence, see step 6 of algorithm 
\ref{algo-svd-gradient}. The parameter \textit{mySeed} allows to have 
reproducible results for different runs on the same matrix.

To show the effectiveness of the proposed SVD method and to illustrate
the choice of the settings, some numerical examples are used.

\subsection{Small size matrix  examples}
Three academic matrix examples,
relations (\ref{eq-example-ab})-(\ref{eq-example-c}),
and the Fisher's iris data,
\cite{Fisher-1936}, were used.
The iris dataset corresponds to the length and the width measurements of
the sepal and the petal for $3$ iris
species: setosa, versicolor and virginica. $50$ measurements were
observed for each species leading to $m=150$ and $n=4$.
\begin{equation}
X_a = \left(\begin{array}{ccc}
  1 &1 &1\\ 0 &2 &1\\ 1 &0 &1\end{array}\right)\mbox{ ; }
X_b = \left(\begin{array}{cccc}
  3 &1 &9 &2\\ 10 &4, &8 &6\\ 7 &6 &12 &1\\ 11 &2 &5 &9\\
  1 &1 &1 &0\end{array}\right)\label{eq-example-ab}
\end{equation}
\begin{equation}
X_c = \left(\begin{array}{ccccc}
22 &10 &2 &3 &7\\ 14 &7 &10 &0 &8\\
-1 &13 &-1 &-11 &3\\ -3 &-2 &13 &-2 &4\\
9 &8 &1 &-2 &4\\ 9 &1 &-7 &5 &-1\\
2 &-6 &6 &5 &1\\ 4 &5 &0 &-2 &2\end{array}\right)\label{eq-example-c}
\end{equation}
The singular values associates with these four examples are shown in Table
\ref{tab-svd-d}. The second and the third eigenvalues for the matrix
$\mathbf{X}_c$ are very close and the last two eigenvalues are zero. 
The matrix $\mathbf{X}_c$ is a challenging example for many SVD
computation algorithms.
\begin{table}[!ht]
\caption{Singular values for the matrices in relations
(\ref{eq-example-ab}) - (\ref{eq-example-c}) and for the iris dataset} 
\label{tab-svd-d}
\begin{center}
\begin{tabular}{rcccc}
\#  &$X_a$ &$X_b$ &$X_c$ &Iris\\ \hline
1 &2.80193774 &26.02508484 &35.32704347 &95.95991387\\
2 &1.44504187 &9.317337969 &20.00000000 &17.76103366\\
3 &0.24697960 &3.298813773 &19.59591794 &3.460930930\\
4 &-                   &0.000000000 &0.000000000 &1.884826306\\
5 &-                   &-                     &0.000000000  &-\\ \hline
\end{tabular}
\end{center}
\end{table}

\begin{table}[!ht]
\caption{Iterations number for small size matrices examples,
the maximum number is set to 200}\label{tab-convergence}
\begin{center}
$q=2$\\
\begin{tabular}{rccccccc}
\hline
$\eta$&0.01&0.1&1  &2 &5 &10 &100\\ \hline
$X_a$ &201 &43 &11 &9 &9 &9 &9\\
$X_b$ &24 &9 &6  &6 &6 &6 &6\\
$X_c$ &179 &179 &179 &179 &179 &179&179\\
Iris  &27 &12 &10 &10 &10 &10 &10\\ \hline 
\end{tabular}
\end{center}
\begin{center}
$q=3$\\
\begin{tabular}{rccccccc} 
\hline
$\eta$ &0.01&0.1 &1  &2 &5 &10 &100\\ \hline
$X_a$  &160&26  &9 &7 &6 &6  &6\\
$X_b$  &8  &5  &5 &5 &5 &5  &5\\
$X_c$  &123&123  &125 &124 &133 &145&201\\
iris   &10 &8  &8 &7 &7 &8 &7\\ \hline
\end{tabular}
\end{center}
\begin{center}
  $q=4$\\
  \begin{tabular}{rccccccc} 
  \hline
  $\eta$ &0.01&0.1 &1  &2  &5 &10 &100\\ \hline
  $X_a$  &85 &16  &7 &6 &5 &7  &7\\
  $X_b$  &5  &4  &4  &4 &4 &4  &4\\
  $X_c$  &140 &178 &201 &201 &201 &201&201\\
  iris   &39 &201 &201 &201 &201 &201 &106\\ \hline
  \end{tabular}
  \end{center}
\end{table}

For a power parameter $q\in \{2,3,4\}$, the step parameter
$\eta$ was varied. The number of iterations required for convergence,
$8$ digits after the decimal point,
are shown in Table \ref{tab-convergence} where $201$ means no
convergence after $200$ iterations. For the examples $X_a$ and $X_b$,
rapid convergence is observed by using $\eta\geq 1$ and $q\in \{2,3,4\}$.
The settings $\eta\leq 0.1$ and $q=3$ seem appropriate for the example 
$X_c$ while $q=3$ and $\eta\geq 1$ are better for the iris dataset.
Using $q=3$ or $q=4$ necessitates 
the same number of operations in the initialization step. 
Using $q=4$ may lead to huge entries in the matrix
$\mathbf{I}_n+\eta(\mathbf{X}^{\mathsf{T}}\mathbf{X})^q$ for some data
examples. Therefore, the choics $q=2,3$ are good compromise. The
choice $q=3$ may be appropriate for large size dataset and
it may be useful to scale the matrix entries before computing the SVD.
The scaling of a nonnegative-valued matrix can consist of dividing each 
entry by the maximum value to have a matrix which entries vary in the 
interval $[0,1]$. The singular values are then divided by the same scalar.

\subsection{High-dimensional example: trachea data}
The trachea dataset corresponds to gene expression profiling using the
sequencing technology. The data was downloaded from the Gene
Expression Omnibus, \cite{Barrett-al-2013}, using the access number
GSE103354,
file: GSE103354\_Trachea\_fullLength\_TPM.txt.gz. 
This dataset is the full-length single cell RNA sequencing for 
$301$ trachea segments of C57BL/6 wild-type
mice, \cite{Montoro-al-2018}. A total of $23,420$ genes (transcripts)
are reported. There are $4,836$ transcripts which have zero count
values for the $301$ cells. These transcripts (rows) were removed and
the remaining sparse dataset has $18,584$ rows, $58\%$
of zero-valued entries, few entries ($0.52\%$) of value greater than 
$1,000$ and a maximum entry value $609,926.3$. 
\begin{figure}[!htb]
  \begin{center}
    \includegraphics[width=0.8\textwidth]{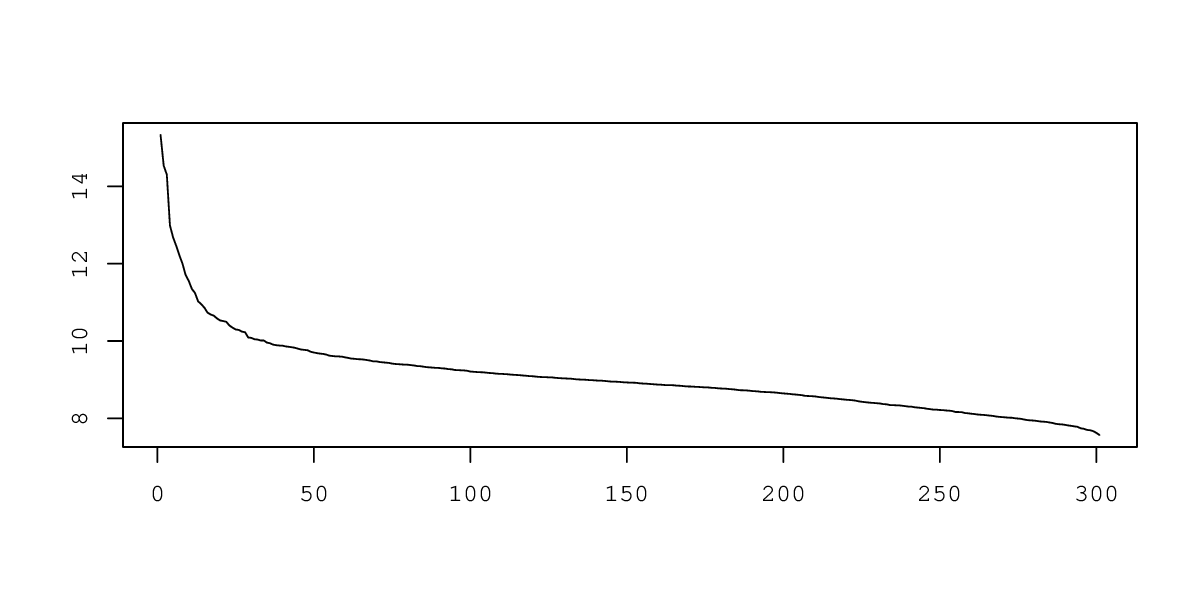} \vspace{-10mm}
    \end{center}
  \caption{Log scale singular values of the trachea dataset}
  \label{fig-eigen-trachea}
\end{figure}

Figure \ref{fig-eigen-trachea} shows a $\log$ scale singular values
behavior for the trachea dataset. 
There is a fast decay for the $20$ first singular 
values and a slow decay for the others. For this example, the
$87$ and $159$ first largest singular values capture $90\%$ and $95\%$ of 
information, respectively.

Five methods available in the literature were used to show the
performance of the proposed method. The calculations were performed on
the same laptop computer equipped with i5-1135G7 processor, 8 GB of
RAM, under Microsoft Windows 11 family and R environment version 4.5.1,
\cite{Rproject}. The \textit{psvd} R package of the
proposed method is available from the comprehensive
R archive network, http://cran.r-project.org/web/packages/psvd/. The
classical SVD method is available through the R build-in \textit{svd()}
function. The Lanczos
algorithm based method for computing all singular values is available
in the \textit{svd} R package, \cite{Korobeynikov-al-2023}.
A C++ library for large scale eigenvalue problem is available in the
Spectra tool, https://github.com/yixuan/spectra. The \textit{RSpectra}
R package allows to use this library for computing the $r$ first 
singular values, \cite{Qiu-al-2024}. The Lanczos bidiagonalization based
\textit{irlba} R package,
\cite{Baglama-al-2022}, also allows to compute the $r$ first singular values.
The probabilistic SVD method is implemented in the \textit{rsvd} R package, 
\cite{Erichson-al-2019}. The $20$, $50$, $100$ or $150$ first largest right 
singular vectors were computed for $5$ methods, see Table 
\ref{tab-svd-d2}. 
The \textit{svd}, \textit{propack.svd}, \textit{RSpectra}, and the \textit{psvd}
methods are almost deterministic and each leads to the same singular
values for different runs on the same matrix. A small error, $\sim$1e-3,
is observed for the \textit{irlba} method singular values for different runs 
on the same matrix, this error is larger, >10, for the \textit{rsvd} method.
These later two methods were run $10$ times, the average results for the
\textit{irlba} method and the best run (minimum mean squares error)
results for the \textit{psvd} method
are reported. Each method is run with the default settings.
For the proposed method, the following settings were used: 
$q$=2, tolerance $\epsilon$=1e-8, $\eta$=1 and \textit{itmax}=200. The
classical \textit{svd()} function was used as reference and three
parameters were examined: the runtime, percentage rate of
reconstruction and the mean squared error (MSE). The MSE value is
obtained by computing the differences between the $r$ first singular values
of a given method and those of the classical SVD method. Then, the mean
of the squares of the differences is computed and shown.
The results obtained are summarized in Table \ref{tab-svd-d2} where the
proposed method appears in bold.
\begin{table}[!ht]
\caption{Comparison of the runtime, the reconstruction rate and the MSE
for six SVD methods}\label{tab-svd-d2}
\begin{center}
\begin{tabular}{crccc}
$r$	&Method      &runtime (s) &reconst.rate \% &MSE\\ \hline
$20$	&svd         &4.12 &80.83 &-\\
	&propack.svd &0.95 &80.83 &$<$1.0e-10\\
	&RSpectra    &0.51 &80.83 &$<$1.0e-10\\
	&irlba       &0.41 &80.83 &2.76e-8\\
	&rsvd        &0.98 &80.83 &889.4\\
	&\textbf{psvd}  &1.61 &80.83 &$<$1.0e-10\\ \hline
$50$	&svd      &4.12 &86.21 &-\\
	&propack.svd &2.81 &86.21 &$<$1.0e-10\\
	&RSpectra    &1.27 &86.21 &$<$1.0e-10\\
	&irlba       &1.34 &86.21 &3.43e-7\\
	&rsvd        &2.14 &86.16 &33,427.6\\
	&\textbf{psvd} &2.51 &86.21 &1.21\\ \hline
$100$ &svd   &4.12 &91.11 &-\\
	&propack.svd &3.75 &91.11 &$<$1.0e-10\\
	&RSpectra    &2.08 &91.11 &$<$1.0e-10\\
	&irlba       &5.42 &91.11 &4.03e-7\\
	&rsvd        &4.51 &91.01 &36,076.5\\
	&\textbf{psvd} &4.37 &91.11 &14.63\\ \hline
$150$ &svd     &4.12 &94.47 &-\\
	&propack.svd &3.94 &94.47 &$<$1.0e-10\\
	&RSpectra    &3.12 &94.47 &$<$1.0e-10\\
	&irlba       &14.22&94.47 &5.3e-3\\
	&rsvd        &7.47 &94.38 &16,203.4\\
	&\textbf{psvd} &6.52 &94.47 &14.35\\ \hline
\end{tabular}
\end{center}
\end{table}

\paragraph{Runtime}
All singular values are simultaneously computed using the classical
\textit{svd} method. When the first $20$ are required,
the \textit{RSpectra}, \textit{irlba}, \textit{rsvd}, \textit{propack.svd} and
\textit{psvd} methods use less runtime compare to the \textit{svd} 
method. For this example, when $r=150$, the
\textit{RSpectra} is the fastest followed by classical \textit{svd},
\textit{propack.svd}, \textbf{\textit{psvd}}, \textit{rsvd} and
\textit{irlba}. The runtime of this later method grows quickly when the
required number of the singular values increases.

\paragraph{Reconstruction rate}
Relation (\ref{eq-svd-reconstr}) was used to compute the
reconstruction rate. The denominator is the same for all methods and
its value is obtained using the classical \textit{svd} method
result. The numerator is computed using the result of each
method. The reconstruction rate is almost the same for all
methods. A very small difference appears for the \textit{rsvd} method
when $r\geq 50$.

\paragraph{Mean squared error}
The differences between the singular values of each method and those
obtained using the classical method where computed and used to 
calculate the MSE. The \textit{propack.svd} and the
\textit{RSpectra} methods have quasi similar reconstruction rate as
the classical \textit{svd} method. 
The first $20$ singular values obtained by all
methods are relatively accurate. The \textit{irlba} method has near
zero MSE value. The \textit{psvd} and \textit{rsvd} methods have more
than zero MSE values. 
The \textit{rsvd} method has relatively higher MSE value and can lead to a 
significant gap between singular values, see Table \ref{tab-svd-d3} 
showing $10$ singular values when $r=100$ (a) or when $r=150$ (b).
\begin{table}[!ht]
\caption{$91$-th to $100$-th singular values obtained using the
classical \textit{svd}, the proposed \textit{psvd} and the
\textit{rsvd} methods: (a) $r=100$, (b) $r=150$}\label{tab-svd-d3}
\begin{center}
  \begin{tabular}{ccccc}
    i &svd &psvd(a) &rsvd(a) &rsvd(b)\\ \hline
    91 &10,847.51 &10,846.65 &10,483.18 &10,808.09\\
    92 &10,804.19 &10,805.04 &10,426.29 &10,754.47\\
    93 &10,638.48 &10,638.33 &10,298.68 &10,603.37\\
    94 &10,600.40 &10,600.56 &10,194.08 &10,540.89\\
    95 &10,413.53 &10,413.28 &10,075.35 &10,371.87\\
    96 &10,376.98 &10,377.18 &9,942.93   &10,329.69\\
    97 &10,332.29 &10,332.33 &9,805.83   &10,264.57\\
    98 &10,289.21 &10,289.21 &9,619.66   &10,225.44\\
    99 &10,195.62 &10,195.62 &9,524.24   &10,137.75\\
    100&10,002.21 &10,002.21 &9,474.41  &9,916.85\\ \hline
    MSE& - &0.161 &231,365 &3,206\\ \hline
  \end{tabular}
\end{center}
\end{table}
Because of the random property of the \textit{rsvd} method, the values
shown in Table \ref{tab-svd-d3} for this method may be different for
different runs while those for the classical method and the proposed methods
remain unchanged.

The \textit{irlba} and the \textit{rsvd} methods should be used for only
large size matrices. These $2$ methods and the RSpectra method are the
fastest for computing a small number of first largest singular vectors.
Only the \textit{propack.svd},
the \textit{RSpectra} and the proposed \textit{psvd} methods can be used 
to compute the $r$ largest singular vectors for any matrix. For all methods 
in the comparison, only the $k$ first singular vectors of the proposed 
method are exactly a neural network autoencoder weight matrix columns. 
Despite its relatively long runtime, the proposed method may be 
useful to have a direct link between data singular vectors and the columns 
of a neural network weight matrix.

\section{Application results}
Truncating of the sum in relation (\ref{eq-svd2}) to $r<n$ 
leads to a $r$ rank matrix $\mathbf{X}$ approximation. The $r$ first largest 
eigenvectors of the matrix $\mathbf{X}^{\mathsf{T}}\mathbf{X}$ are used in 
the principal component analysis and the neural network autoencoder model.

\subsection{Principal component analysis}
In the SVD, no assumption is used on the statistical
properties for the entries of the matrix $\mathbf{X}$.
In practical applications, the entries of the row $\mathbf{x}_i$,
$i=1,2,\ldots,m$, are the measurements associated with 
$n$ attributes or features. Let
us assume that the $n$-dimensional vector $\mathbf{x}$ is a sample 
of a random vector $\mathcal{X}$.
The mean and the covariance of this random vector are:
$\mathbb{E}[\mathcal{X}]\triangleq\boldsymbol{\mu}$ and
$\mathbb{V}[\mathcal{X}] = \mathbb{E}[(\mathcal{X} -
  \boldsymbol{\mu})(\mathcal{X} -
  \boldsymbol{\mu})^T]\triangleq \boldsymbol{\Sigma}$.
The diagonal entries of the matrix $\boldsymbol{\Sigma}$ are the
attributes variances:
$cov(\mathbf{x}_j,\mathbf{x}_j) = \sigma_j^2$, $j=1,2,\ldots,n$. 
The principal component analysis (PCA) consists of the EVD of
the covariance matrix, \cite{Hotelling-1933,Jolliffe-2004}.
When the $\mathcal{X}$'s entries are
heterogenous, there can be a large gap between the amplitudes of the
variances $\sigma_j^2$, $j=1,2,\ldots,n$, leading a misinterpretation 
of the PCA result. To have homogeneous variances,
the samples are standardized and the correlation matrix is
used in the PCA.

\paragraph{PCA for finite number of samples}
When $m$ samples are available, the sample mean $\bar{\mathbf{x}}$ and
the sample variance $s_j^2$ of the $j$-th entry of the sample
$\mathbf{x}_i$, $i=1,2,\ldots,m$, are estimated:
\[\bar{\mathbf{x}} = \frac{1}{m}\sum_{i=1}^m\mathbf{x}_i\mbox{ ; }
  s_j^2 = \frac{1}{m-1}\sum_{i=1}^m(x_{ij} - \bar{x}_j)^2\mbox{ , }
  j=1,2,\ldots, n\]
Then, the entries $x_{ij}^c$ of the centered matrix
$\mathbf{X}_c$ and the entries $x_{ij}^s$ of the standardized matrix
$\mathbf{X}_s$ are obtained as ($i=1,2,\ldots,m$, $j=1,2,\ldots, n$): 
\begin{eqnarray}
  \mathbf{X}_c &:& x^c_{ij} = x_{ij} - \bar{x}_j\nonumber\\
  \mathbf{X}_s &:& x^s_{ij} = (x_{ij} - \bar{x}_j)/s_j\nonumber
\end{eqnarray}
Estimates for the sample covariance and the sample variance matrices are
obtained as $\mathbf{X}_c^{\mathsf{T}}\mathbf{X}_c/m$ and
$\mathbf{X}_s^{\mathsf{T}}\mathbf{X}_s/m$, respectively. Using the
expressions between the SVD factors and those of the EVD,
the PCA result can be obtained by using the matrix
$\mathbf{X}_c$ or the matrix $\mathbf{X}_s$. 

\paragraph{Representation of the samples or the attributes in a reduced space}
The matrices $\mathbf{X}^{\mathsf{T}}\mathbf{X}$ and
$\mathbf{XX}^{\mathsf{T}}$ have the same eigenvalues and related eigenvectors:
\begin{eqnarray}
\mathbf{X}^{\mathsf{T}}\mathbf{Xv} &=& \lambda\mathbf{v}\label{eq-eigen-equiv}\\
\mathbf{XX}^{\mathsf{T}}\mathbf{u} &=& \lambda\mathbf{u}\label{eq-eigen-equiv2}
\end{eqnarray}
where $\mathbf{v}$ and $\mathbf{u}$ are the eigenvectors.
Relation (\ref{eq-eigen-equiv}) allows to write:
$\mathbf{XX}^{\mathsf{T}}\mathbf{Xv} =\lambda\mathbf{Xv}$, meaning that
$\mathbf{u} = \mathbf{Xv}$. Relation (\ref{eq-eigen-equiv2}) can be used in a
similar way to show that $\mathbf{v} = \mathbf{X}^{\mathsf{T}}\mathbf{u}$.
Let $r$ be the reduced dimension used in the PCA, typically $r=2$ or
$3$ for a visual representation. From the PCA result of the matrix
$\mathbf{X}^{\mathsf{T}}\mathbf{X}$, expressions for the samples in the
$r$-dimension attribute space are $\mathbf{X}\mathbf{V}_r$, where
$\mathbf{V}_r$ is a matrix with the $r$ first largest eigenvectors.
Using the expression for the matrix $\mathbf{U}$ in the relation
(\ref{eq-svd-eigen}), the attributes are expressed in the $m$-dimensional
space as:
$\mathbf{X}^{\mathsf{T}}\mathbf{XV}_r\boldsymbol{\Lambda}_r^{-1/2}$, where
$\boldsymbol{\Lambda}_r^{-1/2}$ is the 
square root inverse of the diagonal matrix
associated with the $r$ first largest nonzero eigenvalues.
The PCA result plot for the iris data is shown in Figure
\ref{fig-pca} where the $r=3$ first singular values and singular vectors were
computed, the default settings were used for the other parameters. 
\begin{small}
\begin{figure}[!ht]
\begin{center}
\includegraphics[width=0.95\textwidth]{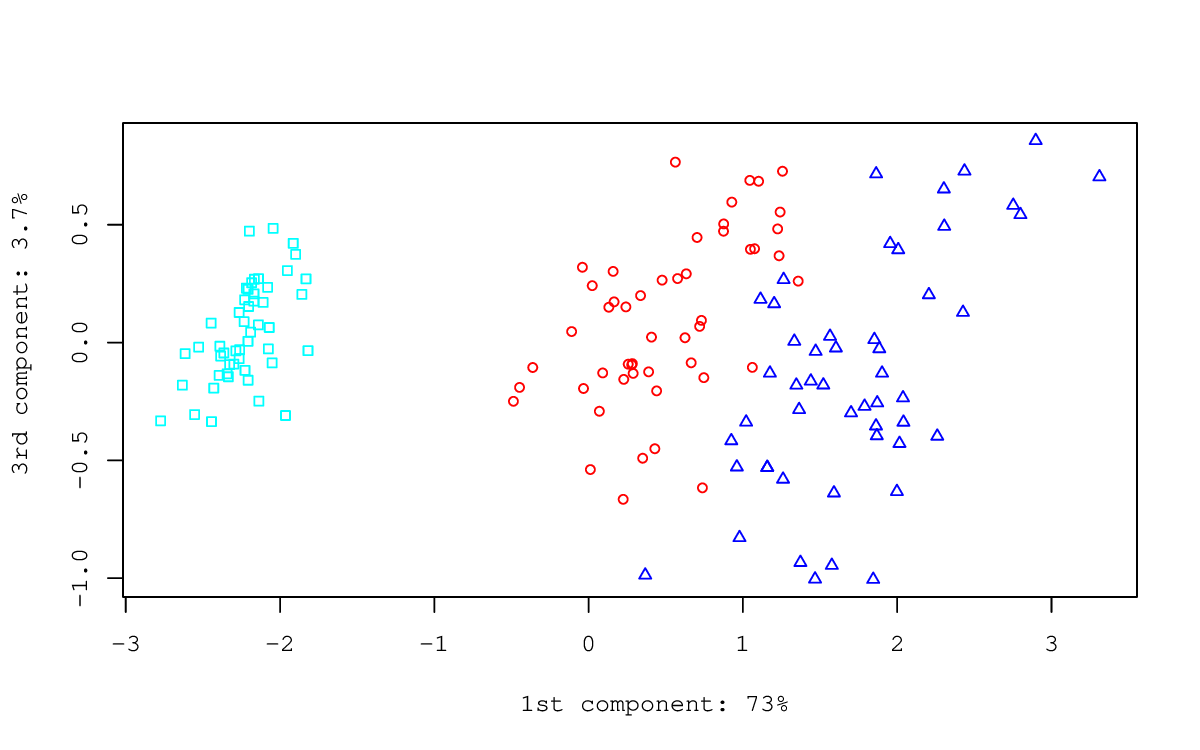} \vspace{-10mm}
\end{center}
\caption{Plot of the iris samples in the first and third
  principal components axes.}
\label{fig-pca}
\end{figure}
\end{small}

\subsection{Autoencoder}
An autoencoder is a neural network (NN) model which allows to
encode input, store the result internally, and finally decode the
stored result to obtain an output. The output is not necessary
strictly equal to the input. In comparison, the SVD leads to the
matrices $\mathbf{U}$, $\mathbf{D}$ and $\mathbf{V}$ which are an
internal representation of the input matrix $\mathbf{X}$. Then, using
the $r$ first largest singular values, a low-rank approximation is obtained:
$\hat{\mathbf{X}} =
\mathbf{U}_r\mathbf{D}_r\mathbf{W}_r^{\mathsf{T}}$.
Thus, the SVD uses linear algebra
operations to solve the dimensionality reduction problem or data
compression. The data can also be compressed
using a feedforward NN model with few layers, 
\cite{Oja-1982,Bourlard-al-1988,Kramer-1991}. A NN
layer is a weighted combination of an input, $\mathbf{x}$,
followed by a transformation function $\varphi(.)$ to obtain an
output: $\mathbf{u} = \mathbf{x}^{\mathsf{T}}\mathbf{W} + \mathbf{w}_o$;
$\mathbf{v} = \varphi(\mathbf{u})$, where $\mathbf{W}$ is a weight
matrix and $\mathbf{w}_o$ is a threshold or bias vector. The
transformation operates on each entry $u$ of the vector
$\mathbf{u}$. Three examples of the transformation function are the
sigmoid, $\sigma(u) = 1/(1+e^{-u})$, the rectified linear unit,
$ReLU(u)=max(0,u)$, and the identity $I(u) = u$. A two layers
feedforward NN model has been used in \cite{Kramer-1991} to perform a
nonlinear SVD for an approximation of $\mathbf{X}$. 
The encoder output $\mathbf{v}_i$ is obtained for each input 
$\mathbf{x}_i$, $i=1,2,\ldots,m$, as: 
\[\mathbf{v}_i = \sigma\left(\mathbf{x}_i^{\mathsf{T}} \mathbf{W}_e +
    \mathbf{w}_{e}\right)\]
where $\mathbf{W}_e$ is the encoder weight matrix of size $n\times r$ and
$\mathbf{w}_e$ is a bias vector of order $r$. The output of the
decoding layer is obtained as: 
\[\hat{\mathbf{x}}_i =
  \sigma\left(\mathbf{v}_i^{\mathsf{T}}\mathbf{W}_n +
    \mathbf{w}_n \right) =
  \sigma\left(\sigma\left(\mathbf{x}_i^{\mathsf{T}}\mathbf{W}_e +
      \mathbf{w}_e\right)^{\mathsf{T}}\mathbf{W}_n+\mathbf{w}_n
  \right)\] 
where $\mathbf{W}_n$ is a weight matrix of size $r\times n$ and
$\mathbf{w}_n$ is a bias vector of order $n$.
The input $\mathbf{x}_i$ and its approximation, the decoder output,
can be used to define a $\ell_2$-norm:
\[\|\mathbf{x}_i - \hat{\mathbf{x}}_i\|\]
From $m$ samples, the
weight matrices, $\mathbf{W}_e$ and $\mathbf{W}_n$, and the bias
vectors, $\mathbf{w}_e$ and $\mathbf{w}_n$ are
obtained after minimization of the
$\ell_2$-norm using a gradient search optimization algorithm.
The approximation matrix obtained using a nonlinear activation feedforward 
NN model is not necessary of reduced rank even the number
of columns of the encoder weight matrix is $r<n$.
However, equivalent approximation to SVD result can be obtained using
the autoencoder model. The difference matrix $\mathbf{X} - 
\hat{\mathbf{X}}$ can be viewed as the noise part of the data and
the approximation $\hat{\mathbf{X}}$ can be used in subsequent
supervised or unsupervised data analysis methods.

\paragraph{Identity transformation function}
The identity transformation function has been used in 
\cite{Plaut-2018,Raff-2022}. The transpose of the matrix $\mathbf{X}$
is used in \cite{Plaut-2018} leading to $\|\mathbf{X} -
\mathbf{WW}^{\mathsf{T}}\mathbf{X}\|$ while the constrained
optimization objective function of the relation (\ref{eq-optim}) is used in
\cite[Chapter 7]{Raff-2022}. In these references, the classical gradient
search optimization method is used to obtain the matrix $\mathbf{W}$.
The $k$ columns of this matrix are not necessary
the $k$ first singular vectors of the data matrix.
For the proposed method, the $k$ first largest singular 
vectors are the $k$ columns of the autoencoder model weight matrix.
When an estimation is available for $\mathbf{W}$, then approximations for 
the matrix $\mathbf{X}$ and the vector $\mathbf{x}_i$ are
$\hat{\mathbf{X}} = \mathbf{X}\hat{\mathbf{W}}\hat{\mathbf{W}}^{\mathsf{T}}$ 
and $\hat{\mathbf{x}}_i = 
\mathbf{x}_i^{\mathsf{T}}\hat{\mathbf{W}}\hat{\mathbf{W}}^{\mathsf{T}}$,
respectively.

\paragraph{MNIST data}
This dataset consists of $60,000$ training images 
and $10,000$ test images, \cite{LeCun-al-1998}. 
Each image of size $28\times 28$ is a sample, a
vector of order $784$. The entries of this vector are pixels which
values are from the set $\{0, 1, \ldots, 255\}$. The behavior of the
singular values and the cumulative sum of the reconstruction
percentage rate of the training images data are show in Figure
\ref{fig-mnist-svd}. The last $72$ singular values are zero. The total
information captured by the $50$, $100$ or $256$ first largest
singular values are $45.41\%$, $60.45\%$ and $82.63\%$, respectively.
\begin{small}
\begin{figure}[!ht]
\begin{minipage}[c]{0.49\linewidth}
\begin{center}
  \includegraphics[width=0.8\textwidth]{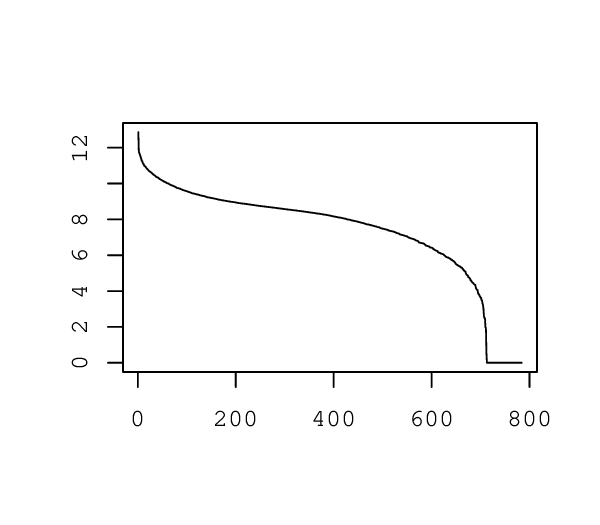}\\ \vspace{-10mm}
  (a) log scale singular values
\end{center}
\end{minipage}
\begin{minipage}[c]{0.49\linewidth}
\begin{center}
  \includegraphics[width=0.8\textwidth]{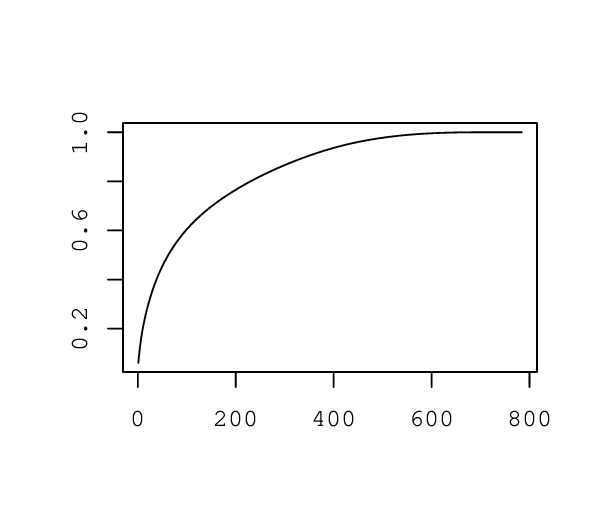}\\ \vspace{-10mm}
  (b) cumulative reconstruction rate
\end{center}
\end{minipage}
\caption{Log scale singular values and cumulative percentage
  reconstruction rate for the MNIST images}
\label{fig-mnist-svd}
\end{figure}
\end{small}

The $256$ first largest singular values were computed using the training
images data and each method. The results are shown in Table \ref{tab-svd-mnist}. 
In the first part of this Table, the default parameter \textit{itmax} (200) is used
while \textit{itmax} = 150(a), 100(b) and 50(c), are used for the 3 last rows
results (proposed method only). Like for the
trachea dataset, the \textit{rsvd} method has the highest MSE value, with a
reconstruction rate almost similar to that of the other methods.
When accurate
singular values are not required, a fast runtime is obtained by using
a small value for the parameter \textit{itmax} in the algorithm
\ref{algo-svd-gradient}. 
\begin{table}[!ht]
\caption{Calculation of the $256$ first singular values for the MNIST images
data using 6 methods: runtime, reconstruction rate and MSE comparison 
results}\label{tab-svd-mnist}
\begin{center}
  \begin{tabular}{rccc}
  method &runtime &reconst.rate &MSE\\ \hline
  svd 	     &85.26 &82.63 &-\\
  propack.svd &71.35 &82.63 &$<$1e-10\\
  RSpectra    &47.50 &82.63 &$<$1e-10\\
  irlba 	     &219.74 &82.63 &5.99e-6\\
  rsvd 	     &87.68 &82.35 &16,631.1\\ 
  psvd 	     &82.75 &82.63 &29.16\\ \hline
  psvd(a)	    &69.72  &82.63 &71.43\\	
  psvd(b)	    &57.71  &82.63 &173.8\\
  psvd(c) 	    &45.35  &82.63 &2,777.4\\ \hline
  \end{tabular}
  \end{center}
\end{table}

The weight matrix obtained using the proposed method
was used to have the $5$ first 
test images approximations. The results are shown in
Figure \ref{fig-mnist-5tests}, where the first column is the
test image and the $3$ following columns are the
results using the $50$, $100$ or $256$ 
first largest singular values approximations, respectively. The
compression rates for these approximations are $15.68$, $7.84$ and
$3.06$, respectively.
\begin{small}
\begin{figure}[!ht]
\begin{center}
\includegraphics[width=0.8\textwidth]{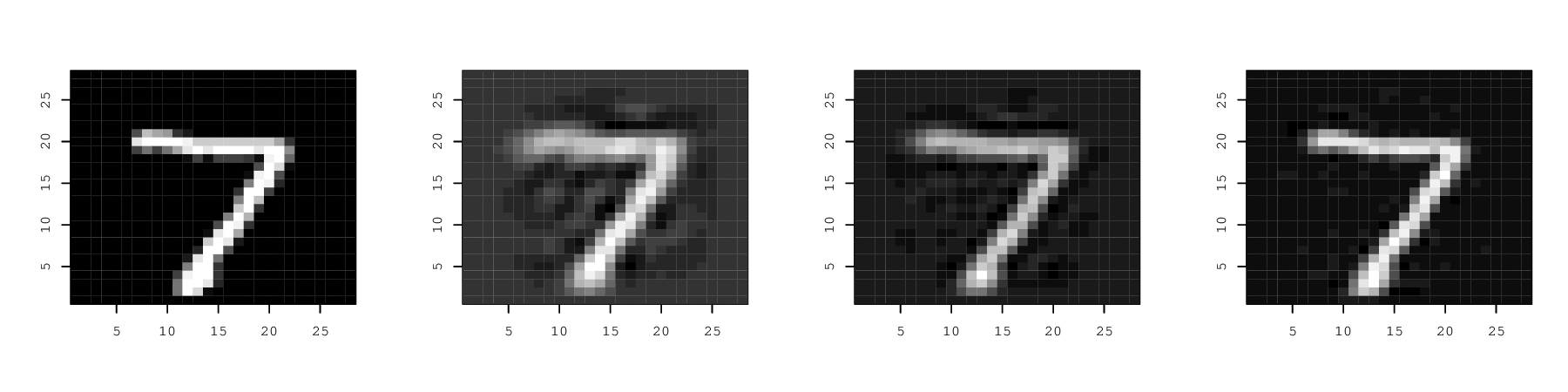}\\ \vspace{-10mm}
\includegraphics[width=0.8\textwidth]{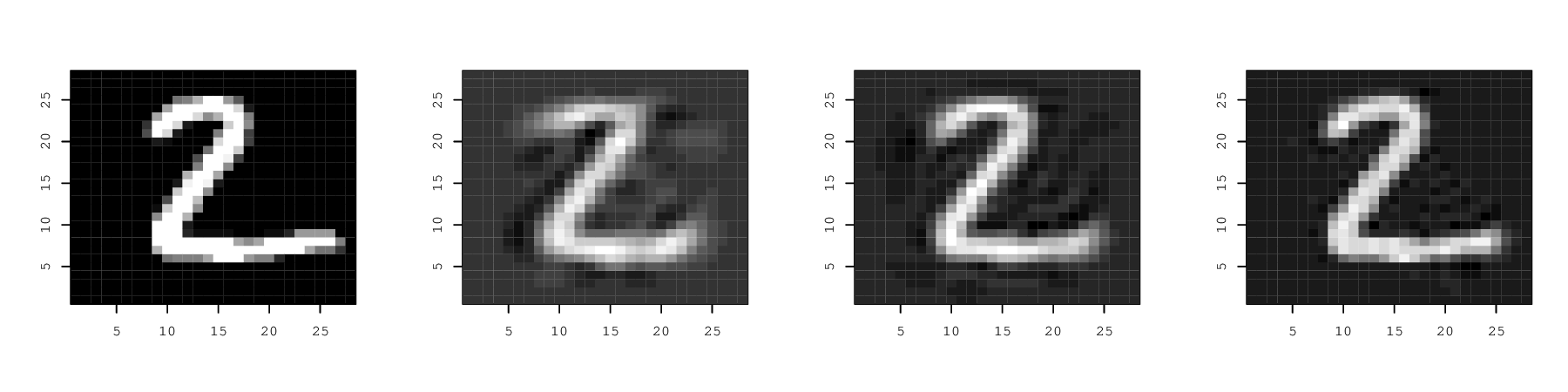}\\ \vspace{-10mm}
\includegraphics[width=0.8\textwidth]{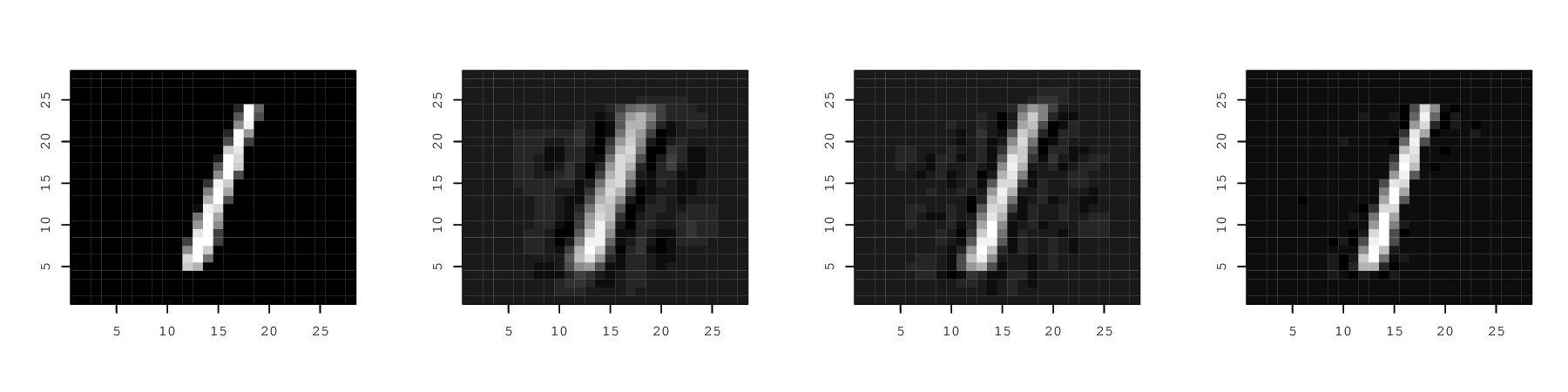}\\ \vspace{-10mm}
\includegraphics[width=0.8\textwidth]{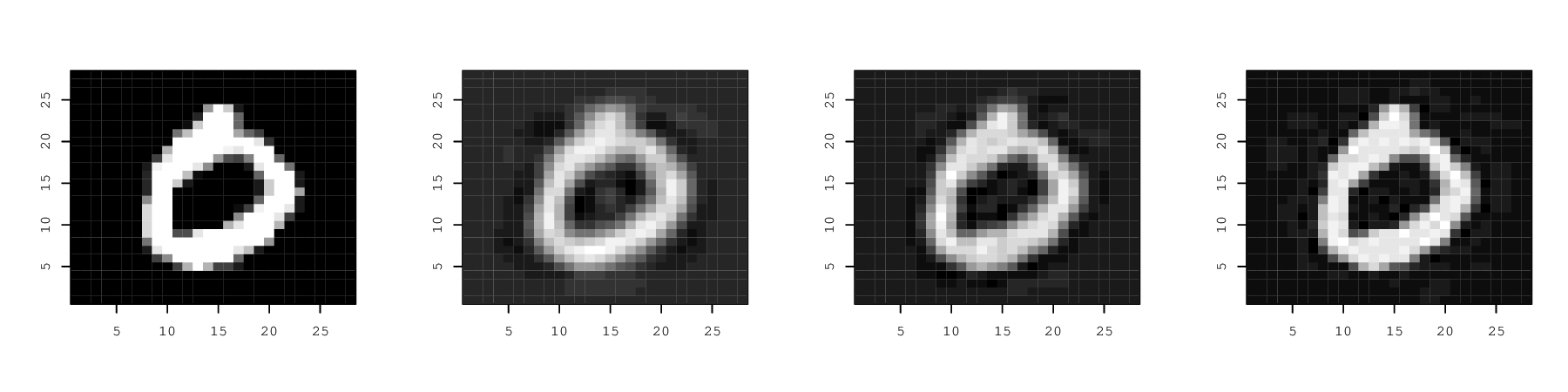}\\ \vspace{-10mm}
\includegraphics[width=0.8\textwidth]{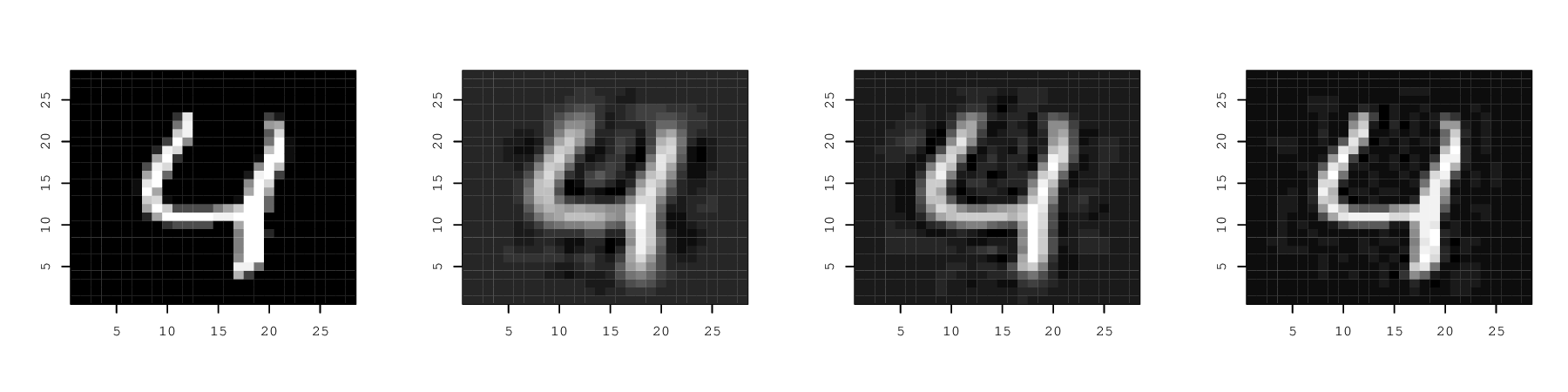}\\ \vspace{-5mm}
\end{center}
\caption{$5$ first test images, from left to right: input, use of the
$50$, $100$ or $256$ first largest singular values to obtain approximation.}
\label{fig-mnist-5tests}
\end{figure}
\end{small}

\section{Conclusion}
A new method is proposed for computing the singular values of a real-valued
entries matrix. This method consists of using a gradient search optimizer 
to minimize a Frobenius norm and allows to compute all or a subset of the $k$ 
first largest singular values and singular vectors. The accuracy of the 
results depends on the user settings. More accuracy may necessitate a long 
runtime depending on the behavior of the singular values decay. 
A parallel implementation may allow to speed up the runtime of the
proposed SVD calculation algorithm. The method 
was applied to principal component analysis and to neural network autoencoder 
model. An independent optimization method is often used in the literature to 
get the singular vectors which are not necessary the $k$ first columns of the
encoder weight matrix unlike for the method proposed.

\appendix
\subsection*{Annexes}
The proposed method is implemented in the \textit{psvd} R package
available on comprehensive R archive network \cite{Rproject}:
http://cran.r-project.org/web/packages/psvd/. 
The power parameter $q$ can be chosen in the set $\{2,3,4\}$, a 
sample from standard normal distribution is used as initial projection 
matrix. The basic power method 
algorithm is given below to show the similarity with the proposed
algorithm. For completeness, a numerically stable Gram-Schmidt
orthogonalization algorithm is provided below. The iterative part of
the proposed algorithm is the more runtime consuming. A C code is used
for this part in the \textit{psvd} R package. Reader can use her/his
preferred computing environment to implement these algorithms.

\subsubsection*{Power method}
\begin{algorithm}
\caption{Power method}
\label{algo-power}
\SetAlgoLined
\KwData{square matrix $\mathbf{A}$}
\KwResult{orthogonal vector $\mathbf{v}$}
Initialization: choose an initial vector $\mathbf{v}^{(0)}$ and compute
  $\mathbf{v}^{(0)}\leftarrow \mathbf{v}^{(0)}/\|\mathbf{v}^{(0)}\|$, set
  tolerance $\epsilon$\;
\For{$t\leftarrow 1,2,\ldots$}{
  Compute a vector $\mathbf{z} = \mathbf{Av}^{(t-1)}$\;
  Update vector: $\mathbf{v}^{(t)}=\mathbf{z}/\|\mathbf{z}\|$\;
  $\delta\leftarrow\|\mathbf{v}^{(t)}-\mathbf{v}^{(t-1)}\|$\;
  \If{$\delta\leq \epsilon$}{stop\;}
}
\end{algorithm}

\subsubsection{Modified Gram-Schmidt algorithm}
\begin{algorithm}[H]
\caption{Modified Gram-Schmidt orthogonalization method}
\label{algo-Gram-Schmidt}
\SetAlgoLined
\KwData{Input = $\mathbf{W}$ has $n$ columns and $m\geq n$ rows}
\KwResult{Output = $\mathbf{V}$, $\mathbf{v}_j$ is the $j$-column of
  $\mathbf{V}$}
$\mathbf{V}\leftarrow \mathbf{W}$\;
\For{$i\leftarrow 1$ \KwTo $n$}{
  $r_i\leftarrow\left\langle\mathbf{v}_i,\mathbf{v}_i\right\rangle^{1/2}$\;
  $\mathbf{v}_i \leftarrow \mathbf{v}_i/r_i$\;
  \If{$i<n$} {
    \For{$j\leftarrow i+1$ \KwTo $n$} {
      $r_j \leftarrow \langle\mathbf{v}_i,\mathbf{v}_j\rangle$\;
      $\mathbf{v}_j\leftarrow \mathbf{v}_j - r_j\mathbf{v}_i$\;
    }
  }
 }
\end{algorithm}
where $\langle\mathbf{x},\mathbf{y}\rangle$ is the inner product of the vectors
$\mathbf{x}$ and $\mathbf{y}$.

\bibliographystyle{abbrv}
\bibliography{psvd3}

\begin{thebibliography}{10}

\bibitem{Arnoldi-1951}
W.~E. Arnoldi.
\newblock The {P}rinciple of {M}inimized {I}terations in the {S}olution of the
  {M}atrix {E}igenvalue {P}roblem.
\newblock {\em Quaterl Appl Math}, 9(1):17--29, 1951.

\bibitem{Baglama-al-2005}
J.~Baglama and L.~Reidel.
\newblock Augmented {I}mplicitly {R}estarted {L}anczos {B}idiagonalization
  {M}ethods.
\newblock {\em SIAM J Sci Comput}, 27(1):19--42, 2005.

\bibitem{Baglama-al-2022}
J.~Baglama, L.~Reidel, and B.~W. Lewis.
\newblock {\em irlba: {F}ast {T}runcated {S}ingular {V}alue {D}ecomposition and
  {P}rincipal {C}omponents {A}nalysis for {L}arge {D}ense and {S}parse
  {M}atrices}.
\newblock Comprehensive R Archive Network, 2022.
\newblock R package, version 2.3.5.

\bibitem{Baldi-al-1989}
P.~Baldi and K.~Hornik.
\newblock Neural {N}etworks and {P}rincipal {C}omponent {A}nalysis: {L}earning
  from {E}xamples {W}ithout {L}ocal {M}inima.
\newblock {\em Neural Net}, 2(1):53--58, 1989.

\bibitem{Barrett-al-2013}
T.~Barrett, S.~Wilhite, P.~Ledoux, and al.
\newblock {NCBI GEO: A}rchive for {F}unctional {G}enomics {D}ata
  {S}ets-{U}pdate.
\newblock {\em Nucleic Acid Res}, 41(Database issue):D991--D995, 2013.

\bibitem{Bishop-2006}
C.~M. Bishop.
\newblock {\em Pattern {R}ecognition and {M}achine {L}earning}.
\newblock Springer, 2006.

\bibitem{Boubekki-al-2021}
A.~Boubekki, M.~Kampffmeyer, U.~Brefeld, and R.~Jenssen.
\newblock Joint {O}ptimization of an {A}utoencoder for {C}lustering and
  {E}mbedding.
\newblock {\em Mach Learn}, 110(7):1901--1937, 2021.

\bibitem{Bourlard-al-1988}
H.~Bourlard and Y.~Kamp.
\newblock Auto-{A}ssociation by {M}ultilayer {P}erceptrons and {S}ingular
  {V}alue {D}ecomposition.
\newblock {\em Biol Cybern}, 59:291--294, 1988.

\bibitem{Brin-al-1998}
S.~Brin and L.~Page.
\newblock The {A}natomy of {L}arge-{S}cale {H}ypertextual {W}eb {S}earch
  {E}ngine.
\newblock In {\em Proc 7th World Wide Web Conf}, pages 107--117, 1998.

\bibitem{Chan-1982}
T.~F. Chan.
\newblock An {I}mprove {A}lgorithm for {C}omputing the {S}ingular {V}alue
  {D}ecomposition.
\newblock {\em ACM T Math Soft}, 8(1):84--88, 1982.

\bibitem{Dembele-2021}
D.~Demb\'el\'e.
\newblock A {M}ethod for {C}omputing the {P}erron {R}oot for {P}rimitive
  {M}atrices.
\newblock {\em Numer Linear Algebra Appl}, 28(1):e2340, 2021.

\bibitem{Demmel-al-1999}
J.~Demmel, M.~Gu, S.~Eisenstat, and al.
\newblock Computing the {S}ingular {V}alue {D}ecomposition with {H}igh
  {R}elative {A}ccurary.
\newblock {\em Linear Algebra Appl}, 299:21--80, 1999.

\bibitem{Demmel-al-1990}
J.~Demmel and W.~Kahan.
\newblock Accurate {S}ingular {V}alues of {B}idiagonal {M}atrices.
\newblock {\em SIAM J Sci Stat Comput}, 11(5):873--912, 1990.

\bibitem{Eckart-al-1936}
C.~Eckart and G.~Young.
\newblock The {A}ppromixation of {O}ne {M}atrix by {A}nother of {L}ower {R}ank.
\newblock {\em Psychometrika}, 1(3):211--218, 1936.

\bibitem{Erichson-al-2019}
N.~B. Erichson, S.~Voronin, S.~L. Brunton, and J.~N. Kutz.
\newblock Randomized {M}atrix {D}ecompositions {U}sing {R}.
\newblock {\em J Stat Soft}, 89(11):1--48, 2019.

\bibitem{Fan-al-2024}
Z.~Fan, C.~Yang, B.~Lin, Y.~Yang, and Q.~Shi.
\newblock Convex {S}et-{O}riented {S}ingular {V}alue {D}ecomposition with
  {B}ounded {U}ncertainties.
\newblock {\em J Comput Appl Math}, 448:115942, 2024.

\bibitem{Fard-al-2020}
M.~M. Fard, T.~Thonet, and E.~Gaussier.
\newblock Deep k-{M}eans: {J}ointly {C}lustering with k-{M}eans and {L}earning
  {R}epresentations.
\newblock {\em Pattern Recognit Lett}, 138:185--192, 2020.

\bibitem{Fisher-1936}
R.~A. Fisher.
\newblock The {U}se of {M}ultiple {M}easurements in {T}axonomic {P}roblems.
\newblock {\em Ann Eugenics}, 7:179--188, 1936.

\bibitem{Frieze-al-2004}
A.~Frieze, R.~Kannan, and S.~Vempala.
\newblock Fast {M}onte {C}arlo {A}lgorithms for {F}inding {L}ow-{R}ank
  {A}pproximations.
\newblock {\em J ACM}, 51(6):1024--1041, 2004.

\bibitem{Frobenius-1912}
F.~G. Frobenius.
\newblock Ueber {M}atrizen aus nicht {N}egativen {E}lementen.
\newblock {\em Sitzungsberichte der K\"{o}niglich Preussischen Akademie der
  Wissenschaften}, Berlin:456--477, 1912.

\bibitem{Golub-al-1965}
G.~H. Golub and W.~Kahan.
\newblock Calculating the {S}ingular {V}alues and {P}seudo-{I}nverse of a
  {M}atrix.
\newblock {\em SIAM J Num Anal}, 2(2):205--224, 1965.

\bibitem{Golub-al-1970}
G.~H. Golub and C.~Reinsch.
\newblock Singular {V}alue {D}ecomposition and {L}east {S}quares {S}olutions.
\newblock {\em Numer Math}, 14:403--420, 1970.

\bibitem{Golub-al-2000}
G.~H. Golub and H.~A. {van der V}orst.
\newblock Eigenvalues {C}omputation in the 20th {C}entury.
\newblock {\em J Comput Appl Math}, 123:35--65, 2000.

\bibitem{Golub-al-1996}
G.~H. Golub and C.~F. {van L}oan.
\newblock {\em Matrix {C}omputations}.
\newblock The Johns Hopkins Univ Press, Baltimore, 3rd edition, 1996.

\bibitem{Gu-2015}
M.~Gu.
\newblock Subspace {I}teration {R}andomization and {S}ingular {V}alue
  {P}roblems.
\newblock {\em SIAM J Sci Comput}, 37(3):1139--1173, 2015.

\bibitem{Halko-al-2011}
N.~Halko, P.~G. Martinsson, and J.~A. Tropp.
\newblock Finding {S}tructure with {R}andomness: {P}robabilistic {A}lgorithms
  for {C}onstructing {A}pproxilate {M}atrix {D}ecompositions.
\newblock {\em SIAM Rev}, 53(2):217--288, 2011.

\bibitem{Hanson-al-1969}
R.~J. Hanson and C.~L. Lawson.
\newblock Extension and {A}pplications of the {H}ouseholder {A}lgorithm for
  {S}olving {L}inear {L}east {S}quares {P}roblems.
\newblock {\em Math Comput}, 23(108):787--812, 1969.

\bibitem{Hestenes-al-1951}
M.~R. Hestenes and W.~Karush.
\newblock A {M}ethod of {G}radient for the {C}alculation of the
  {C}haracteristic {R}oots and {V}ectors of a {R}eal {S}ymmetrix {M}atrix.
\newblock {\em J Res Natl Bur Stand}, 47(1):45--61, 1951.

\bibitem{Horn-al-2019}
R.~R. Horn and C.~R. Johnson.
\newblock {\em Matrix {A}nalysis}.
\newblock Cambridge Univ Press, 2nd edition, 2019.

\bibitem{Hotelling-1933}
H.~Hotelling.
\newblock Analysis of a {C}omplex of {S}tatistical {V}ariables into {P}rincipal
  {C}omponents.
\newblock {\em J Educ Psychol}, 24(6):417--441, 1933.

\bibitem{Jolliffe-2004}
I.~T. Jolliffe.
\newblock {\em Principal {C}omponent {A}nalysis}.
\newblock Springer, 2nd edition, 2004.

\bibitem{Kampffmeyer-al-2019}
M.~Kampffmeyer, S.~L{\o}kse, F.~M. Bianchi, L.~Livi, A.-B. Salberg, and
  R.~Jenssen.
\newblock Deep {D}ivergence-{B}ased {A}pproach to {C}lustering.
\newblock {\em Neural Net}, 113:91--101, 2019.

\bibitem{Knyazev-al-1994}
A.~V. Knyazev and A.~L. Skorokhodov.
\newblock Preconditioned {G}radien-{T}ype {I}terative {M}ethods in a {S}ubspace
  for {P}artial {G}eneralized {E}igenvalue {P}roblem.
\newblock {\em SIAM J Numer Anal}, 31(4):1226--1239, 1994.

\bibitem{Korobeynikov-al-2023}
A.~Korobeynikov and R.~M. Larsen.
\newblock {\em svd: {I}nterfaces to {V}arious {S}tate-of-{A}rt {SVD} and
  {E}igensolvers {S}oftware}.
\newblock Comprehensive R Archive Network, 2023.
\newblock R package, version 0.5.7.

\bibitem{Kramer-1991}
M.~A. Kramer.
\newblock Nonlinear {P}rincipal {C}omponent {A}nalysis {U}sing
  {A}utoassociative {N}eural {N}etworks.
\newblock {\em AIChE J}, 37(2):233--243, 1991.

\bibitem{Lanczos-1950}
C.~Lanczos.
\newblock An {I}terative {M}ethod for the {S}olution of the {E}igenvalue
  {P}roblem of {L}inear {D}ifferential and {I}ntegral {O}perators.
\newblock {\em J Res Nat Bur Standards}, 45(4):255--282, 1950.

\bibitem{Langville-al-2005}
A.~N. Langville and C.~D. Meyer.
\newblock A {S}urvey of {E}igenvector {M}ethods for {W}eb {I}nformation
  {R}etrival.
\newblock {\em SIAM Rev}, 47(1):135--161, 2005.

\bibitem{Larsen-1998}
R.~M. Larsen.
\newblock Lanczos {B}idiagonalization with {P}artial {R}eorthogonalization.
\newblock {\em DAIMI Report Series}, 27(537):1--101, 1998.

\bibitem{LeCun-al-2015}
Y.~{Le C}un, Y.~Bengio, and G.~Hinton.
\newblock Deep {L}earning.
\newblock {\em Nature}, 521(7553):436--444, 2015.

\bibitem{LeCun-al-1998}
Y.~{Le C}un, L.~Bottou, Y.~Bengio, and P.~Haffner.
\newblock Gradient-{B}ased {L}earning {A}pplied to {D}ocument {R}ecognition.
\newblock {\em Proc IEEE}, 86(11):2278--2324, November 1998.

\bibitem{Li-al-2025}
X.~Li, X.~Wu, T.~Wang, Y.~Xie, and F.~Chu.
\newblock Fault {D}iagnosis {M}ethod for {I}mbalanced {D}ata {B}ased on
  {A}daptive {D}iffusion {M}odels and {G}enerative {A}dversarial {N}etworks.
\newblock {\em Eng Appl Artif Intell}, 147:110410, 2025.

\bibitem{Liberty-al-2007}
E.~Liberty, F.~Woolfe, P.-G. Martinsson, V.~Rokhlin, and M.~Tygert.
\newblock Randomized {A}lgorithms for the {L}ow-{R}ank {A}pproximation of
  {M}atrices.
\newblock {\em Proc Natl Acad Sci, USA}, 104(51):20167--20172, 2007.

\bibitem{Lu-al-2022}
H.~Lu, C.~Chen, H.~Wei, Z.~Ma, K.~Jiang, and Y.~Wang.
\newblock Improved {D}eep {C}onvolutional {E}mbedded {C}lustering with
  re-{S}electable {S}ample {T}raining.
\newblock {\em Pattern Recognit}, 127:108611, 2022.

\bibitem{Marcus-al-1964}
M.~Marcus and H.~Minc.
\newblock {\em A {S}urvey of {M}atrix {T}heory and {M}atrix {I}nequalities}.
\newblock Allyn and Bacon, Inc., Boston, 1964.

\bibitem{Meyer-2000}
C.~D. Meyer.
\newblock {\em Matrix {A}nalysis and {A}pplied {L}inear {A}lgebra}.
\newblock SIAM, Philadelphia, 2000.

\bibitem{Montoro-al-2018}
D.~T. Montoro, A.~L. Haber, M.~Biton, and al.
\newblock A {R}evised {A}irway {E}pithelial {H}ierarchy {I}ncludes
  {CFTR-E}xpressing {I}onocytes.
\newblock {\em Nature}, 560(7718):319--324, 2018.

\bibitem{Mrabah-al-2020}
N.~Mrabah, N.~M. Khan, R.~Ksantini, and Z.~Lachir.
\newblock Deep {C}lustering with a {D}ynamic {A}utoencoder: {F}rom
  {R}econstruction {T}owards {C}entroids {C}onstruction.
\newblock {\em Neural Net}, 130:206--228, 2020.

\bibitem{Ogita-al-2020}
T.~Ogita and K.~Aishima.
\newblock Iterative {R}efinement for {S}ingular {V}alue {D}ecomposition {B}ased
  on {M}atrix {M}ultiplication.
\newblock {\em J Comput Appl Math}, 369:112512, 2020.

\bibitem{Oja-1982}
E.~Oja.
\newblock A {S}implified {N}euron {M}odel as {P}rincipal {C}omponent
  {A}nalyzer.
\newblock {\em J Math Biol}, 15:267--273, 1982.

\bibitem{Perron-1907}
O.~Perron.
\newblock Zur {T}heorie der {M}atrices.
\newblock {\em Mathematiche Annalen}, 64(2):248--263, 1907.

\bibitem{Petersen-al-2012}
K.~B. Petersen and M.~S. Petersen.
\newblock The {M}atrix {C}ookbook.
\newblock http://matrixcookbook.com, November 15, 2012.

\bibitem{Plaut-2018}
E.~Plaut.
\newblock From {P}rincipal {S}ubspaces to {P}rincipal {C}omponents with
  {L}inear {A}utoencoders.
\newblock {\em arXiv}, page 1804.10253, 2018.

\bibitem{Qiu-al-2024}
Y.~Qiu, J.~Mei, G.~Gunnebaud, and J.~Niesen.
\newblock {\em {RS}pectra: {S}olvers for {L}arge {S}cale {E}igenvalue and {SVD
  P}roblems}.
\newblock Comprehensive R Archive Network, 2024.
\newblock R package, version 0.16-2.

\bibitem{Rproject}
{R Core Team}.
\newblock {\em R: {A L}anguage and {E}nvironment for {S}tatistical
  {C}omputing}.
\newblock R Foundation for Statistical Computing, Vienna, Austria, 2024.

\bibitem{Raff-2022}
E.~Raff.
\newblock {\em Inside {D}eep {L}earning: {M}ath, {A}lgorithms {M}odels}.
\newblock Manning, Shelter Island, NY, USA, 2022.

\bibitem{Saad-1980}
Y.~Saad.
\newblock Variations on {A}rnoldi's {M}ethod for {C}omputing {E}igenelement of
  {L}arge {U}nsymmetric {M}atrices.
\newblock {\em Linear Algebra Appl}, 34:269--295, 1980.

\bibitem{Sleijpen-al-1996}
G.~L.~G. Sleijpen, A.~G.~L. Booten, D.~R. Fokkema, and H.~A. van~der Vorst.
\newblock Jacobi-{D}avidson {T}ype {M}ethods for {G}eneralized {E}igenproblem
  and {P}olynomial {E}igenproblems.
\newblock {\em BIT Numer Math}, 36(3):595--633, 1996.

\bibitem{Stewart-1993}
G.~W. Stewart.
\newblock On the {E}arly {H}istory of the {S}ingular {V}alue {D}ecomposition.
\newblock {\em SIAM Review}, 35(4):551--566, 1993.

\bibitem{Tropp-al-2018}
J.~A. Tropp, A.~Yurtsever, M.~Udell, and V.~Cevher.
\newblock Practical {S}ketching {A}lgorithms for {L}ow-{R}ank {M}atrix
  {A}pproximation.
\newblock {\em arXiv}, pages 1609--00048, 2018.

\bibitem{Trosten-al-2024}
D.~J. Trosten, S.~Lokse, R.~Jenssen, and M.~Kampffmeyer.
\newblock Leveraging {T}ensor {K}ernels to {R}educe {O}bjective {F}unction
  {M}ismatch in {D}eep {C}lustering.
\newblock {\em Pattern Recognit}, 149:e110229, 2024.

\bibitem{Xie-al-2016}
J.~Xie, R.~Girshick, and A.~Farhadi.
\newblock Unsupervised {D}eep {E}mbedding for {C}lustering {A}nalysis.
\newblock In {\em Proc 33rd ICML}, page 10p. JMLR W\&CP, volume 48, 2016.

\bibitem{Xu-1993}
L.~Xu.
\newblock Least {M}ean {S}quare {E}rror {R}econstruction for
  {S}elf-{O}rganizing {N}eural-{N}ets.
\newblock {\em Neural Nets}, 6(5):627--648, 1993.

\bibitem{Yang-al-2024}
C.~Yang and Q.~Shi.
\newblock An {I}nterval {P}erturbation {M}ethod for {S}ingular {V}alue
  {D}ecomposition ({SVD}) with {U}nknown-{B}ut-{B}ounded ({UBB}) {P}arameters.
\newblock {\em J Comput Appl Math}, 436:115436, 2024.

\bibitem{Zhao-al-2024}
D.~Zhao, W.~Cai, and L.~Cui.
\newblock Adaptive {T}hresholding and {C}oordinate {A}ttention-{B}ased
  {T}ree-{I}nspired {N}etwork for {A}ero-{E}ngine {B}earing {M}onitoring
  {U}nder {S}trong {N}oise.
\newblock {\em Adv Eng Inf}, 61:102559, 2024.

\bibitem{Zhou-al-2007}
Y.~Zhou and Y.~Saad.
\newblock A {C}hebyshev-{D}avidson {A}lgorithm for {L}arge {S}ymmetric
  {E}igenproblem.
\newblock {\em SIAM J Matr Anal Appl}, 29(3):954--971, 2007.

\end{thebibliography}
\end{document}